\documentclass[a4paper,10pt]{article}
\usepackage[utf8]{inputenc}
\usepackage{amsmath,amssymb}
\newcommand{\vsp}{\vspace{.5cm}}

\begin{document}

\title{Some remarks on monogenity}
\author{Tove Dahn}
\maketitle
\begin{abstract}
 We will discuss representations of monogenic functions over very regular groups. In particular, related to normal models. Further, we discuss codimension, maximality and monodromy for the representations.
\end{abstract}
\subsection{Summable distributions}
\label{sec:summable}
We will discuss several representations, that depend of choice of topology.
Consider for instance $f \in \mathcal{D}_{L^{1}}'$ (\cite{Schwartz66}), implies $\mathcal{F}f=Q f_{0}(x,y)$ with $f_{0}$ continuous in a nbhd (neighborhood) of singularities.
In particular, when singularities are on on the diagonal, $f_{0}$ can be chosen of negative type, outside the diagonal. When further Q is reduced, we have that f corresponds to a very regular action. We also consider the Laplace transform $\mathcal{L}(T)=\mathcal{F}(e^{-<\xi,x>} T) = Q' f_{0}$, with $Q',f_{0}$ as above. (\cite{Schwartz52})

Given $\log f \in L^{1}(L)$, with algebraic singularities on complex lines L, starting with a very regular domain (\cite{BrelotChoquet51}), there is a $f \in \mathcal{D}_{L^{1}}'(L)$, with $Q f_{0}=\mathcal{L}(\log f)$, that is $\log f \in L^{1}(G_{2})$ is approximated by  $\log f_{0} \in H(G_{8})$ (\cite{Dahn24}). Note that algebraic representations exclude spirals. f=0 on a set of positive measure iff $f_{0}=0$, on the same set. Given $f_{0}=e^{\phi}$, with $\phi >0$ implies $f_{0} \neq 0$, in particular f has algebraic zeroes on a pseudo convex domain defined by $\{ \phi >0 \}$. Every $T \in H'(\mathbf{R}^{2})$ has representation in $\mathcal{E}^{'(0)}$ over $H(G_{8})$ and when $\phi \in H$, $\mathcal{L}T(\phi)=T(e^{\phi})$ (\cite{Lelong68}).

\vsp

Assume $T \in \mathcal{D}_{L^{1}}'$, note that the pseudo topology, corresponds to completion of
the space of test functions to a strict condition, that is $\mathcal{D}_{L^{1}} \rightarrow \dot{B}$. Simultaneously, the polar set gives completion to canonical topology.
In particular, given the restriction R to a Hausdorff normed space (injective), then ${}^{t}R$ relative canonical dual space is surjective and is used to define the polar set.
Example: $dV \in G$ and
$VT(V\phi)=UT(\phi) \rightarrow V \phi$ continuous , for instance $P \rightarrow 1/P$ absolute continuous, with $VT \sim 1/P$ for a polynomial P (geometric equivalence), that is VT has measure zero singularities and satisfies a strict condition, when $1/P \rightarrow 0$.
Example: $T-I \in C^{\infty}$, then ``spectrum'' is given by traces (diagonal), characterized by conjugation.

Concerning non-canonical topology: consider $<T(\phi),dU>=<T,dU>(\phi)$, for $T \in X'$. Thus, the condition $U X' \subset X'$ is necessary, or ${}^{t}U X \subset X$. Completion to a strict condition, according to $\mathcal{D}_{L^{1}}' \subset (\mathcal{D}_{L^{1}})'$, can be given by ${}^{t} U \phi \in \dot{B}$, when $\phi \in \mathcal{D}_{L^{1}}$ and $T \in (\mathcal{D}_{L^{1}})'$,
with $UT \in \mathcal{D}_{L^{1}}'$. The polar set can be given by $\phi \notin L^{1}$, but $\phi \in L^{1}(dU)$, that is $dI \notin G_{2}$, a discontinuous group.

\vsp

In the plane, when $C^{\infty}(G_{2}) \simeq H(G_{8})$, we have representations in $\mathcal{E}^{'(0)}$ (convolution algebra). Since the order of zero, is determined by the order for the restriction to lines, it is under a strict condition, sufficient  to consider a pluri-complex representation. The order of the representation is determined by contraction (section \ref{sec:contraction}).
Note that completion to $L^{1}(G_{2})$, with a strict condition implies it is sufficient to consider a representation over $H(G_{2})$.
$\mathcal{L}T \in $ Exp, determines a linear functional  $T \in H'$. When $e^{\phi} \in H$, $T(e^{\phi})$ gives a linear representation.
Thus, $\phi \in H(G_{8})$, with $\mathcal{L}T(\phi)=T(e^{\phi})$, can be represented for $e^{\phi} \in H(G_{8})$.
Example: continuation $\mathcal{S} \rightarrow \dot{B}$, with ${}^{t} U \phi \in \dot{B}$, gives a linear representation  $ UT \in \mathcal{D}_{L^{1}}'$ (pseudo topology). Thus, we have $\widehat{UT}=P(\xi) f_{0}$, where $f_{0}$ is continuous and locally bounded outside traces (the diagonal).

\newtheorem{lemma1}{Lemma 1}[section]
\begin{lemma1}
Assume $f=e^{\psi} \in \mathcal{H}$=Exp $L^{1}$, where $\mathcal{H}$ is a dual topology to a space between  $\mathcal{D}_{L^{1}}$ and $L^{1}$. Then we have existence of $\phi$ in this space, such that $<\psi,dV>,<\phi, d {}^{t} V>$ are linear, for $d V \in G_{\mathcal{H}}$.
\end{lemma1}

Given V algebraic, $\phi$ can be chosen as dual to $f \in \mathcal{H}$.
Assume $d V=\rho dI$ with $ \rho_{1} \psi \in L^{1}$ (transpose of $\rho$) and consider $< \psi, {}^{t} \mathcal{L}I(\phi)>$, with $\psi$ subharmonic and completed to its best harmonic majorant (\cite{AhlforsSario60} ch.II), for instance
$\rho_{1} \tilde{\psi}$  in $L^{1}$, corresponding to restriction of $\phi$ to a strict condition. Thus, $\mathcal{L}I \in H_{0}'$ ($H_{0}$ are harmonic functions with removable singularities). Note that the Fourier transform of $e^{-<\xi,x>} e^{\psi}$ gives the Laplace transform, $< \psi, {}^{t} \mathcal{L}I(\phi)>$, over for instance $\phi \in \mathcal{S}$ (\cite{Schwartz52}).

\vsp

When $T \in H'$, we can write $\mathcal{L}(T)(\phi)=T(e^{\phi})$ over $\phi \in H$ (cf \cite{Lelong68}). Example: $\phi(u,v) \in H(G_{2})$ does not imply $\phi(x,y) \in H$, but $\phi \in C(G_{2})$ implies $\phi \in H(G_{8})$ in the plane.
Assume $S(e^{\tilde{\phi}})=T(\tilde{\phi})$, with $\mathcal{L}(S)=T$. We assume monodromy in the sense that given $\phi \in H$,
when $\phi(u,v) \in H$, we have local existence of (u,v), single valued between two elements in the range to G, that is between elements accessible for regular approximations. This assumes convexity in the representation space.

Assume $UF \in H'$ and consider $V \phi$, completion to $H_{0}$. Sufficient for ${}^{t}UV=I$ $(= V {}^{t}U)$ over $H_{0}$, is that $d I \in G_{H_{0}}$, that is ${}^{t} U =V^{-1}$ over $H_{0}$.
Over $H(G_{8})$ a sufficient condition to determine ${}^{t}U$ uniquely, is that $dI \in (G_{8})_{H}$ has point support, that is dI is extremal.
Example: given $d U \in G$ and conjugation according to $UT(e^{V f})=T(e^{f})$, for some V with $dV \in G_{2}$, alternatively
$UT(e^{Vf})=WT(e^{f})$, with $WT \in H'$.
Example: when $T(e^{\phi} - I)=\widehat{T}(\phi) - \widehat{T}(0)$, given $e^{\phi} - I \in H$, we have that $T \in H'$.

\subsection{Nuclearity and continuous groups}
\label{sec:nuclear}

Starting with Dirichlet-series $f(s)=\Sigma a_{n} e^{-\lambda_{n}s}$, with $s=\sigma + it$, in particular ordinary Dirichlet-series, where $\lambda_{n}=\log n$, we discuss the problem when $1 > \mid a_{n} \mid > \epsilon$, for some terms, but $\Sigma \mid a_{n} \mid^{\mu} < \infty$, for some $\mu$. According to  (\cite{Riesz15}, sec 2.5), when $\overline{A}_{n}=\Sigma_{1}^{n} \mid a_{j} \mid$ and $S=\Sigma \mid a_{n} \mid e^{-\lambda_{n}s}$, we assume existence of $\overline{\sigma}=\lim \frac{\log \overline{A}(n)}{\lambda_{n}}$ ($\leq \mu$). Then S is absolute convergent only for $\sigma > \overline{\sigma}$. In particular, given S convergent in $s=s_{0}$, S has the limit $f(s_{0})$ along every path for which $\mid am(s-s_{0}) \mid \leq \alpha < \pi/2$ (\cite{Riesz15}, Theorem 5). Further, $(-1)^{\rho} f^{(\rho)}(s) \sim \Sigma a_{n} \lambda_{n}^{\rho} e^{-\lambda_{n} s}$, uniformly convergent, for a positive entire number $\rho$ and $\sigma \geq \sigma_{0}$, for some $\sigma_{0}$ (\cite{Riesz15}, Theorem 4).

When the results are applied to the representation space (u,v), we have that the convergence domain determines a neighborhood of invariants. The convergence domain can then be used to determine type of movement. When the domain is symmetric and when the conjugation gives a projective decomposition of the plane, we have convergence on a disk.

On a convergence domain, we have $\lim_{s \rightarrow s_{0}}f(s)=f(s_{0})$. Assume $R(G)=\Omega$ defines the domain for $f \in C^{0}(\Omega)$. When $G_{2} \subset G$ has $d I \notin G_{2}$, $G_{2}$ contains a trajectory outside the domain for convergence. In the plane for s, we can assume $H(G_{8})$ dense in $C^{0}(G_{2})$. Starting with a normal model (\cite{Dahn23}), $d U \rightarrow (dU,dV) \rightarrow dN(U,V)$, we must have $dN/dV \rightarrow 0$ (or $dV/dU \rightarrow 0$) regularly at the boundary. Simultaneously $dN/dU \rightarrow 1$, defines a domain for convergence in $C^{0}(G_{2})$. Thus,  when $s_{0} \in \mbox{ reg } \Omega$ connected (points regular for H), using $G_{2} \subset G_{8}$, we have a regular approximation property. When the continuation  is irreducible, we have that dN(U,V)=dI implies dU=dI.

Assume C a closed convex cone in $\mathbf{R}^{n}$. C can be defined  by the convex closure of finitely many  extremal points (Minkowsky). An extremal ray $ S \subset C$,
according to $S=\mathbf{R}_{+} v$ for $v \neq 0$, is such that $\forall x,y \in C$ such that  $x+y \in S$, we have that $x,y \in S$.

A contraction domain refers to, (df/dy)(dy/dx)=df/dx regular (single valued).
 Assume $F(\gamma)(x,y,z)=f(x,y,z)$ and $\gamma(u,v) \in H$, a regular approximation property relative (dU,dV), gives that the pseudo base $\gamma$ is regular in (x,y,z). A very regular distribution E, means that $E(\gamma)$ defines a regular $C^{\infty}$ neighborhood of evaluation.

Given G continuous, (u,v) open implies that $C^{\infty}_{c}(u,v)$  decomposable.
Example: consider $\int f(v)dv=\int f(v(u))dv$. Assume $v=\varphi(u)$, with $f \circ \varphi$ regular, when f(u) regular.
Given $f=f_{1}(u) \otimes f_{2}(v)$ and $\sigma=dv/du$
$\int \int f(u,v) du dv = \int f_{1} d u \int f_{2} d v$
$=\int f_{1}(u) du \int \varphi^{*} f_{2} \sigma d u$ . On the other hand, if $(u,v) \rightarrow u-v$ is dense (u,v regularly situated), we have $< \delta(u-v),f(u,v)> = \int f(u,u) du$.

\vsp
For a TVS (topological vector space) X, an irreducible component is a maximal subspace.  Example: assume $X \subset Y$ with Y TVS and X connected, presence of a polar set then implies X is considered as irreducible in Y.
Example: $X=X(G_{2})$ with $dU+dV=dI$ and X maximal in Y.
A continuous image of for instance an irreducible boundary, is irreducible.

Assume $J_{j}$ a proper ideal defined by $A_{j}-\lambda I$, for a linear operator $A_{j}$ in $\mathcal{D}_{L^{1}}'$. $J_{j}$ is maximal iff $A_{j}$ are irreducible.
Example: $A_{1}A_{2}=I$, for instance $A_{1}(dU)$,$A_{2}(dV)$, is interpreted such that  $A_{1}A_{2}$ surjective and every two elements can be connected using U,V, in particular $G_{2}$ is continuous. Example:
Assume $dU+dV=dI$, that is a projective decomposition is regarded as irreducible, when invariants are given by dU=dI or dV=dI. Spirals are interpreted as joint spectrum.

\subsection{Monogenity}
\label{sec:monogenic}
Denote with C the class of monogenic functions (\cite{Borel12}).
Multivalentness  for $f'$ generates several (say two) classes of normals $n_{1},n_{2}$, given $(n_{j})$ generates a vector space, we have $f \notin C$. Example: when bounded variation implies determined tangent (with respect to arc length), this implies determined normal, but not necessarily monogenity.
Further, $C \subset C^{\infty}$ implies $\mathcal{E}' \subset C'$. We note that multivalentness  does not exclude a representation in $C'$. Example: $d U=\rho d I$, with $\rho$ absolute continuous, defines an absolute continuous measure when $\rho=1$.
Using Sard's theorem (\cite{Malgrange66}) $f \in C^{\infty}$
on an open set $\subset \mathbf{R}^{n}$, implies $D^{k}f=0$ has measure zero for some $k \geq 1$, that is we have abscense of non-trivial first surfaces.

\vsp

 When $f'$ is single valued on a trajectory $\gamma$, this does not imply $(\gamma,f(z))$  single valued. Assume f=P/Q, with Q polynomial. Consider $\phi : x \rightarrow y$ absolute continuous , this preserves algebraicity. When f has representation according to Weierstrass (\cite{Borel12}), the same holds for $f \circ \phi$. Monogenic functions extend the class of Weierstrass-functions, in particular under the condition $\phi(u)=v$ reversible, for instance through a normal model. Example: semi-algebraic sets, when Q reduced,  $\{ P < \lambda \} \subset \{ P/Q < \lambda' \}$, that is we have semi-algebraic subsets to $\{ f < \lambda \}$ and a reduction to a monogenic model.
Example: every polynomial Q has restriction to a reduced domain $\Omega_{1}$.
When $\phi :\Omega_{2} \rightarrow \Omega_{1}$ is absolute continuous , it preserves algebraicity. In particular invertibility on a reduced domain. Example: $\phi(x)=1/x$ absolute continuous , then $Q \circ \phi$ is polynomial in $\infty$, given Q polynomial in 0, cf. preserves a constant value.

\vsp

Given R(dU) defines a domain of holomorphy, f can be continued analytically for inner points, where bd R(dU)=$\{ d U=dI \}$. Further, for instance when f(u)=f(0)=0 at bd R(dU), assume dV=0 at the boundary, that is $d V \neq dI$, when $f(0,v) \neq 0$. When V,U are assumed of different type, then V defines a continuation of the domain to U.
In this manner, f is continuable by every regular movement with dV=0 on the boundary and monogenity can be of higher order (cf. section \ref{sec:involution} and resolution of sng).
 Assume instead $f \in C(G_{2})$, then $H(G_{8})$ defines possibly multivalent continuation.
 Monogenic continuation implies that R(dU) can be continued to a planar disk by dV,
  that is using sequential continuation. Presence of traces implies abscense of monogenity.
  Example:  Assume dU+dV=dI
 and existence of a trajectory $\eta$, with $dU \rightarrow (dU,dV)$ a continuation along $d V(\eta)=0$. When V absolute continuous, $\eta$ is invariant to V and U simultaneously, why $\eta=0$ or $\{U-I=0 \} \cap \{ dU-dI =0 \} = \emptyset$.

 $UFU^{-1}=F$ corresponds to $U \rightarrow {}^{t}U$ does not change type. Example:
 $U \mid f \mid^{2}=Uf U g$, where $g=\overline{f}$. Example: $\mid f \mid=1$ implies $U(1/f)=U^{-1}(f)$
 (cf. section \ref{sec:exponential} approximative representations), that is U(f)U(1/f)=1.

 \subsubsection{Reduction to the boundary}
Assume $d U \rightarrow (dU,dV) \rightarrow dN$ defines analytic continuation through a normal model. Necessary for a normal model, is that dN/dV=0 (or dV/dU=0) at the boundary, in particular given dN algebraic in dU, we have monodromy and monogenity.

\newtheorem{lemma2}[lemma1]{Lemma 2}
\begin{lemma2}
Assume $f(v_{1},v_{2}) \in C^{\infty}$, with $v_{1} \rightarrow v_{2}$ absolute continuous and linearly  independent. Then we have $\lim_{v_{2} \rightarrow 0} f(v_{1},v_{2})= f(v_{1})$, that is given that the domain is defined by the boundary (convexity) $v_{2}=0$, we can consider $f(v_{1})$ as monogenous
\end{lemma2}

Assume $d v_{2}/ d v_{1} \mid_{v_{2}=0}=0$ and $f_{v_{1}} + f_{v_{2}} dv_{2}/ d v_{1}=f_{v_{1}}$ in $v_{2}=0$.

$f(u,v) \rightarrow f(u,v/u) \rightarrow f(u,1)$, that is density  outside $u/v=1$ implies existence of measures $\mu(u,v)=0$ on u=v.
v/u=1 is trivially rectifiable in the plane, for instance with respect to arc length $s=\int_{0}^{1} \sqrt{1 + (dv/du)^{2}} d t$ when v=v(u). Note that contact transforms preserve traces (zero-lines).
Example: $T'(f)=\mid f'(x) \mid$ continuous and single valued. Then given bounded variation, T(f) is monogenous.

\vsp

On a parabolic surface, we have for $f=e^{\phi}$, that $\phi < 0$ implies $\phi=const$ (\cite{AhlforsSario60}).
Assume $\tilde{\phi}$ harmonic $=e^{\psi}$, with $\log \log \mid f \mid \leq C \mid x \mid$, for instance $\psi \in L^{1}$ linear, this implies $\phi$ has algebraic singularities. Example:
When P polynomial, from Laplace $\mathcal{L}I(P)=e^{P}$. Assume $\int \mid d F \mid < \infty$ on a planar disk D, with $dF=e^{P} d x$ and where $\Gamma=\{ e^{P}=const \}$ is rectifiable, that is zero sets in bd D of measure zero are mapped onto sets of measure zero in $\Gamma$ under $dF \rightarrow e^{P}$.

\newtheorem{lemma3}[lemma1]{Lemma 3}
\begin{lemma3}
Assume $\phi$ absolute continuous , and f monogenic, then we have that $f \circ \phi$ monogenic.
 \end{lemma3}

Assume $\phi$ absolute continuous . Thus of bounded variation, that is $T(\phi)<\infty$, in particular the derivative is finite. $d \phi/dx=0$ implies $\phi=const$, implies $(f \circ \phi)'=0$. Let $\phi=\phi_{+}-\phi_{-}$, with $\phi_{\pm}$ non-decreasing (\cite{Riesz56}).  When $d \phi_{\pm}>0$, $\phi_{\pm}$ is monotonous and thus  single valued and
$\frac{df}{d \phi_{\pm}} \frac{d \phi_{\pm}}{dx}$ is single valued, that is $f'$  single valued implies $(f \circ \phi)'$  single valued.
Monogenity implies that df has a unique definition on $\Gamma \rightarrow \phi(\Gamma)$.
In particular, f(x,y) monogenic, has continuation through a normal. Example: $F(\gamma)(\zeta)=f(\zeta)$ with $\gamma$ absolute continuous , that is f monogenous in $\zeta$, implies $F(\gamma)$ monogenous in $\zeta$.
Example: $d N \bot dU \times dV$, where
$u \neq 0$ and $v \neq 0$, that is defines a surface and when further $u \neq v$, we have a unique normal.
Example: $d N_{1} \bot \Gamma_{1}$ and $d N_{2} \bot \Gamma_{2}$. Assume $\phi(\Gamma_{1})=\Gamma_{2}$, with $\phi$ absolute continuous, then $d N_{1}$ algebraic implies $d N_{2}$ algebraic.

\subsubsection{Partial monogenity}
$\int f d U=\int f d I$ does not imply dU=dI, that is $\Sigma=\{ x \quad Uf(x)=f(x) \}$ and $d \Sigma=\{ x \quad d U(x)=dI(x) \}$ can be different. Consider $f(u,v/u) \rightarrow f(u,0)$, when $v/u \rightarrow 0$. Further $f(u,v) = f_{1}(u) \otimes f_{2}(v)$, for instance in $C^{\infty}((u,v))$, for (u,v) open, that is $f(u,v/u)=\lim f(u,0) \otimes f(0,v/u)$, where  decomposability assumes monogenity.

We consider partial monogenity $f(u,v,w)=g(u,v)h(w)$.
Necessary for this is decomposability for $((u,v),w)=(z,w)$, in particular when f is partially monogenic, $f \in C^{\infty}(z,w)$ does not imply $f \in C^{\infty}(z)$. Example: consider (z,w) and F=sng f, with
$W F \cap F = \emptyset$, but $(U,V) F \cap F \neq \emptyset$, that is an one-parameter regularization.

Starting from (u,v), consider contraction to u. Given linear contraction ($L^{1}$ in the phase) according to a strict condition, the polar set is removable. Starting from $C^{\infty}(G_{2}) \simeq H(G_{8})$, where $G_{2}=(dU,dV)$, given for instance $dW \notin G_{2}$ and $\frac{\delta f}{\delta w} \frac{\delta w}{\delta u} + \frac{\delta f}{\delta w} \frac{\delta w}{\delta v} \neq 0$ and regular. Then, we have reduction to $H(G_{2})$, that is the polar set is a removable set on R(dW).
Example: v does not affect the integral in dudv, when v(u)=const. For linear continuation $(u,v_{1}) \rightarrow (u,v_{2})$, we assume $d v_{1}/d v_{2}$ homotopical with constant mapping, so that traces are removable. Linear contraction is a continuous (rectifiable) deformation, contraction with traces is discontinuous .
Starting from a normal model, according to UN=NU and $dN(U,V) \rightarrow dI$ at the boundary and given f continuous, real over a removable segment, we have that df(u,v) $\sim$ d f(u,n).
Example: $<U^{\bot} f,d V>=<{}^{t} V U^{\bot} f,dI>=<f,d V \times dU^{\bot}>$. Assume the limit when $d V \rightarrow dI$ is = the limit when $d U^{\bot} \rightarrow d I$, where invariants are assumed to have a regular (pseudo convex) neighborhood, that is given that $R(d U)^{\bot}$ is generated by $d(U^{\bot},V)$, it is necessary for a regular normal to $R(dU,dU^{\bot})$, that $\delta f / \delta v \neq 0$.
Consider $d U^{\bot} \rightarrow dV$ with common invariants. Note that $\frac{\delta f}{\delta v} \neq 0$ does not imply $\frac{\delta f}{\delta u^{\bot}} \neq 0$ (diffuse front), in particular monogenity in dV does not imply monogenity in $d U^{\bot}$.

 \vsp

$d (U + V)^{2} = 2 (U+V) d (U+V)$, that is dUV=U dV + V d U.  In particular given $R(UdV) \subset R(dV)$ and $R(VdU) \subset R(dU)$, we have $R(d UV) \subset R(d U + dV)$.
Assume $UdVf=UVg$ for some  g and $VdUf=VUh$, for some  h. Example: $\int f d UV + \int_{\Gamma}df(u,v) = \frac{1}{2} UVg + \frac{1}{2} VUh$, note that when d UV=d VU, that is when df symmetric in (u,v), we can choose g,h in the same ideal, for instance g(u,v)=h(v,u).
Thus, it is necessary for dUV(f)=0, that $g \rightarrow h$ continuous through involution.
Given $d U \rightarrow dV$ harmonic conjugation, analyticity is preserved in the plane by a pure mapping. Note that monogenity is preserved by a pure mapping, since $f^{\diamondsuit}=-if$ (harmonic conjugate).

\subsection{Spectral theory}
\label{sec:spectral}
Assume $f \in H(G \times G)$ and consider transformations (V,U)f on a Banach space $\widehat{X}$, for instance $\tilde{f}=f(\phi) \in X$, where $\phi \in \mathcal{D}$.
When we assume that $\tilde{f}(u,v) \in C^{\infty}$ implies $\tilde{f}(u) \in C^{\infty}$, then $sng(\tilde{f}(u)) \subset sng(\tilde{f}(u,v))$. Consider $\tilde{f}$ locally, as a constant coefficients differential operator in X. Given a spectral resolution of dI, from $\{ \tilde{f}(u,v) < \lambda \}$ in $\widehat{X}$, the spectrum to the operator defined by $\tilde{f}(\xi)$, can be determined. An algebraic action of v implies $\sigma(\tilde{f}_{\lambda}(u,v)(D))=\sigma(\tilde{f}_{\lambda}(u)(D))$.

Example $f=e^{\phi}$ with $\phi'$
linear =L(x,y). When  $\phi_{x}=L(x,\cdot)$ and in the same manner for $\phi_{y}$, we have $\phi=\alpha <x,x> + \beta <y,y> + \gamma$, for constants $\alpha,\beta,\gamma$. On a triangulated 2-dimensional surface (\cite{AhlforsSario60}), we can represent the domain to $\phi$ using $G \times G$, with a projective decomposition.

Example: When $U \mathcal{L}I(\phi)=\mathcal{L}I(V \phi)$, given V determines the type for U uniquely, we have when $\phi \in H$, that $\mathcal{L}I(V \phi) \rightarrow dU$ locally injective, here the duality is according to Fourier. Assume a pseudo convex domain, $e^{\phi}$ with $\phi>0$ close to the boundary.
Assume the domain in (u,v), where v is determined by an involution condition and u is determined by duality (conjugation ). Then $f \in C^{\infty}(u)$ is dependent of type of movement, that is dependent of whether U preserves sng f or not.

  Consider continuation without affecting the spectrum, for instance algebraic continuation.
  A quasi-inverse is defined by common invariant sets $\Sigma$, for $dU,dU^{-1}$. This means that $\Sigma(dUV)=\Sigma(dV)=\Sigma(dI)$, that is does not define disk-neighborhoods. On Hausdorff uniformities (\cite{Bourbaki89}), $\Sigma(dUV)$ can be used to define joint spectrum, when $dU,dV$ linearly independent.

Constant surfaces are regarded as singularities for $C^{\infty}$.
Assume $\{ f=c \} \subset$ a domain of holomorphy.
When $\{ f=c \}$ an analytic set and $\{ f=c \}^{c}$ can be regarded locally as an analytic set, then nbhd $\{ f=c \}$ is
Borel (Suslin).  On a transversal L to f=c, sng are regarded as isolated points (Sard). When f=c corresponds to $\delta_{x}(f)=c$, supp $\delta_{x}$ is regarded as analytic (holomorphic leaves).
Assume $\{ df=0 \} \subset \{ \mid df \mid =0 \}$, that is T(f)=const. defines a neighborhood of first surfaces.
Example: $\Sigma=\{ dU=dI \}$ and $\Gamma=\{ f=c \}$, assume $\Sigma \rightarrow \Gamma$ continuous. Then sng can be described using G, that is starting from  sng $\subset \Gamma$,  $\Sigma \rightarrow \Gamma$ approximates sng.
Example $\Delta T(f)=0$ defines a harmonic measure, thus a parabolic boundary. The corresponding domain can be approximated by very regular f.

Assume $\Omega^{c}$ in a domain of holomorphy, where $\Omega$ analytical.
Assume $\phi \in H(\Omega)$ and algebraic locally on $W \subset \Omega^{c}$, then there is an extension to W. When W is analytic, it can be regarded as removable, for instance a zero line to an analytic measure.

\newtheorem{lemma5}[lemma1]{Lemma}
\begin{lemma5}
 When $d_{\Gamma} T(f) <C$ defines a neighborhood of $\Gamma=\{ d_{\Gamma}=0 \}$, where $d^{2}_{\Gamma}$ algebraic, then $\Gamma$ can be regarded as removable.
\end{lemma5}

Example: $d_{\Gamma} \mid f \mid < C$, given f polynomial, defines a domain of holomorphy. Assume f of bounded variation on $\Omega$, for an open set $\Omega$, where $T(f)$ defines the distance to the ideal boundary and $d_{\Gamma}$ denotes the distance to $\Gamma=\{ f=const \}$, with $d^{2}_{\Gamma}$ algebraic, that is $d^{2}_{\Gamma} \sim p \geq 0$. Then $T(f) < C/d_{\Gamma} \sim C d_{\Gamma}/p$, when $d_{\Gamma}$ close to 0 and f is rectifiable outside $p=0$. Further, when g has bounded variation on ($\Omega - \Gamma$) and $gh \in C^{\infty}(\Omega-\Gamma)$, where h is flat on $\Gamma$, we have existence of $\tilde{gh} \in C^{\infty}(\Omega)$, that is $d_{\Gamma}=0$ (here $d_{\Gamma}=h$) defines $\Gamma$, which here is removable (\cite{Malgrange66}).

\subsubsection{Joint invariants}
When $\mathcal{H} \ni f \rightarrow c$ in $\infty$, we have that the
spectral function to constant surfaces
corresponding to $f = c$, has cluster sets.
Example: $f(\phi)$ regularization, where $f(\phi) \equiv c$, when $\phi$ close to $\delta$, that is does not have an approximations property in $\mathcal{H}$, that is $d I \notin G_{\mathcal{H}}$. According to the lifting principle, when $\phi(x,y)$ polynomial, $\mathcal{H} \ni f(\phi) \rightarrow \phi \in H$ is continuous.
By extending G, we can consider $\mathcal{H} \ni f(\phi) \rightarrow \phi \in H(G_{8})$.
In particular, $f(\phi) \in \mathcal{H}$, with $\mathcal{H} \subset L^{1}$ under a strict condition, motivates the reduction $G_{8} \rightarrow G_{2}$ continuous.
Example: $f_{v} v_{u} \neq 0$ regular, where $dV /d U =\rho \in L^{1}$, $\rho \rightarrow 0$ in $\infty$, that is where $f_{v} \neq 0$ is regular and $\rho \in L^{1}$, we have a contraction domain, according to $G_{8} \rightarrow G_{2}$ continuous.
Simultaneously,  when f is independent  of v in some direction, we do not have a contraction domain.
When f=c on $\Omega$, the surface is invariant for every $(u,v) \Omega \subset \Omega$.

Assume $f \in H(K(G_{4}))$ (\cite{Dahn24}) in a neighborhood of $f(u,v)=const$, then there is not contraction to $G_{2}$ in H.
More precisely, when $(dU^{1},d V^{1}) \in K(G_{4})$, we assume $(dU^{1},dV^{1})=(\rho_{1} dU,\rho_{2} dV)$, with $\rho_{j} \neq const$, $j=1,2$.
f(u,v)=f(0,0) defines joint invariants on a disk neighborhood,
that is (0,0) is not isolated in a disk neighborhood.
$C^{\infty}_{c}(\Omega_{1} \times \Omega_{2})$ is  decomposable (\cite{Treves67}), when $\Omega_{j}$ are open sets in $\mathbf{R}^{n}$, . Example: $\{ (0,0) \} \subset \subset D$, a disk neighborhood
(cf. Hausdorff). Thus, $C^{\infty}_{c}$ is not decomposable over D, given joint invariants.
When $u \rightarrow v$ continuous and f(u,v)=f(0,0) implies u=v, we have that Hausdorff (\cite{Bourbaki89}) implies $C^{\infty}_{c}$ decomposable outside u=v.

Example: assume $\Sigma$ is defined by $f(u,v)=c$ iff $f(u,0)$ and $f(0,v)$ are constant.
When $f_{u} \neq 0$, for u close to 0, we are outside $\Sigma$. Further, when $\phi(u)=v$ is absolute continuous, with $\phi_{u}=0$, as $v \rightarrow 0$, we have $f_{u}(u,v)=f_{u}(u,0)$.
Assume dUV=dI, where dU,dV of different type, that is defines a disk neighborhood.
$dI \notin G$ with $dUV=dI$ does not imply dU=dI. Assume conversely $dUV=dI$ implies dU=dI, then also dV=dI, and an irreducible boundary implies joint invariants. Example: dI(u)d I(v)=dI(u,v) according to f(0,0)=f(u,v) defines joint invariants.
Example: $Uf=\int f d U + \int_{\Gamma} df$ and $VUf=\int \lambda f dI + \int_{\Gamma_{\lambda}} df$.
When $\Gamma_{\lambda}$ is a zero line to df, we have joint invariants, that is defined by (u,v).

\subsection{Removable set}
\label{sec:removable}
Assume $\phi : v_{1} \rightarrow v_{2}$ absolute continuous , such that  $d \phi=0$ implies
$\phi=const$, that is does not affect the action, $\frac{df}{d \phi} \frac{d \phi}{d v_{1}}=0$, when $d \phi=0$, that is $d f$ independent of $v_{2}=const$.
In particular, where $d \phi=0$, we have $\lim_{v_{1} \rightarrow 0} d f(v_{1})= \lim_{v_{1} \rightarrow 0} d f(v_{1},v_{2})$, that is monogenity. Thus, $v_{1} \rightarrow v_{2}$ absolute continuous, represents an equivalence, independent  of initial value.

Define $\tilde{\Gamma_{2}}$ by $\Gamma_{1} \rightarrow \Gamma_{2}$ absolute continuous, where $\Gamma_{1}$ satisfies a strict condition. Consider $\Gamma_{1} \subset \Gamma_{1} \cup \tilde{\Gamma}_{2} \subset \Gamma_{1} \cup \Gamma_{2}$. Example: $U^{\bot} f = \int f d U^{\bot} + \int_{\tilde{\Gamma}} df=f$, given $d U^{\bot}=0$ implies $U^{\bot}=I$, for instance $f \in L_{ac}^{1}(\tilde{\Gamma})$

\vsp

Example: $\frac{d}{d x} \log f=f'(x)/f(x)=\rho$ regular, implies $f(x)=(1/\rho)(x)f'(x)$.
$f'/f=const$ corresponds to multivalentness for $\frac{d}{d x} \log f$.
When $\int f d I=\int f'(x) d x$, monogenity implies that f is single valued relative dI or weakly single valued.
Consider g(x)=f(x,y), where y=y(x) is linear, then $g'(x) \sim (1+y')f_{x}$, that is $g'(x)dx \sim (dx + dy)f_{x}$, thus $\int_{\Gamma} df=\int_{\Gamma} g_{x} dx$. Consider $\Gamma=\Gamma_{x} \times \Gamma_{y}$, for instance $\Gamma$ open, that is given $f_{x}$ symmetric with respect to $(x,y) \rightarrow (-y,x)$, we have $\int_{\Gamma} df=0$ is not dependent on dy, that is y can be regarded as a parameter.

Assume $d U \rightarrow d V$ a single valued continuation, such that  (dU,dV) gives a strictly larger domain.
$d I \notin R(dU)$, when $dI \in R(dU,dV)$, implies presence of the corresponding proper sub-ideal R(dU), that is abscense of an approximation property for R(dU). Assume f not regular on dU=dI.
When $dV=\rho dU$ and $\rho \in L^{1}$, on an unbounded domain, we must have $\rho \neq const$.
Example: for a monogenic f, when (dU,dV) in a very regular subgroup of $G \times G$ (\cite{Dahn24}), is a single valued continuation of dU, we must have $\rho$ regular, $\rho' \neq 0$, .

\subsubsection{Polar subgroup}

Consider the polar set in G. Assume $G_{2} \subset G=G_{2} \bigoplus H$, with $d I \notin G_{2}$.
Then we can regard $H$ as envelop of normals in G. Given $D^{k}f(h_{1},h_{2},\ldots)=0$ on a set of measure zero, when k sufficiently  large (Sard), we have density for $G_{2}$.
Example: given $f=f_{1}(u) \otimes f_{2}(h)$, where $f_{j}$ are analytic, then zero-sets in $G_{2}$ can be regarded as Borel. Example: $H \rightarrow G_{2}$ absolute continuous, preserves sets of measure zero.

Assume $G=G_{2} \bigoplus H$ and $f=f_{1}(u) \otimes f_{2}(h)=0$ as above. When $(u,h)$ irreducible, we have $f_{1}(0,h)=0$ or
$f_{2}(u,0)=0$, that is when f symmetric, an analytic set. This assumes density for the representation, in particular $f=0$ for some  $(u,v)$. Given $u \rightarrow v$ quasi-inverses, that is does not change type, we have that (u,v) is not dense in the plane.
(\cite{Julia19}).

\vsp

Assume $(u,v) \bot n$ and $(u,0) \bot n$. Further, given ((u,v),n)=(u,(v,n)), then we have $u \bot (v,n)$, when $v \rightarrow 0$. cf. Weyl: in $L^{2}$, $n \bot u$,$n \bot v$ implies $n \in C^{1}$. Assume $n \bot \Omega$, where $\Omega$ is
symmetric in u,v and includes traces.
Assume $n \bot u$, an unique normal, further $u \bot v$ and $n \bot v$, then $n \bot (u,v)$.

Assume $(u_{1},u_{2},u_{3})$ linearly independent and $(u_{1},u_{2})$ conjugated according to Fourier, with $u_{3} \bot (u_{1},u_{2})$, as $u_{2} \rightarrow 0$ and $u_{1}$ represents translation. Assume $n \bot u_{3}$, as $u_{2} \rightarrow 0$, transversal according to a normal tube. Then n corresponds to translational movement parallel to $u_{1}$.

\vsp

Example: assume $\Gamma$ =hyperboloid. Given $\Gamma$ is represented by G, the dimension for $\Gamma$ is dependent of the order of G. When G is given by translation, $\Gamma$ is 0-dimensional, when G is given by translation, rotation, $\Gamma$ is 1-dimensional (\cite{Dahn24}). Assume $X_{+},X_{-}$ on either side of $\Gamma$. They are regularly situated, if
$ 0 \rightarrow \mathcal{E}(X_{+} \cup X_{-}) \rightarrow_{\pi_{1}} \mathcal{E}(X_{+}) \bigoplus \mathcal{E}(X_{-}) \rightarrow_{\pi_{2}} \mathcal{E}(X_{-} \cap X_{+}) \rightarrow 0$ is exact (\cite{Malgrange66}). Sharp fronts implies  decomposability for open sets. Assume $f(\Gamma)=0$, for $f=g_{1}(X_{-}) - g_{2}(X_{+})$, cf. very regular representation. When $\Gamma$ is algebraic, f can be chosen as polynomial.

\vsp

Example:  $d U \rightarrow d U \times dV$ is projective on traces, but dU-d V=0 does not imply dU=0.
Conformal mappings map geodetics on geodetics.
Given $\gamma$ geodetic, there is an $\eta$, with $\Delta \eta=0$ and $\phi^{-1}(\eta)=\gamma$, for a mapping $\phi$ (\cite{Brelot49}). Example: max is assumed on the boundary $\Gamma$, but the boundary does not define the domain.
Assume $\gamma$ a line in the representations space, $\phi$ a continuous mapping $\gamma \rightarrow $ in 3-space. Then $\eta$ can be a spiral
(\cite{Dahn19} $d U_{S}^{\diamondsuit}$), that is $\eta$ (involutive image of) has a harmonic definition, but $\eta$ does not separate the 3-space regularly.

\vsp

The diagonal $\Delta$ is defined by supp dI. Assume
$\mu_{t} - \delta \in C^{\infty}$, with $\mu_{t}$ a harmonic measure, then $\mu_{t}=0$ defines a parabolic boundary.
$T \in \mathcal{D}^{' F}$ very regular implies $T-I \in C^{\infty}$, this implies existence of a TVS X, such that  $T \in C^{\infty}(X \backslash \Delta)$, where $\Delta$ denotes the diagonal.
Given $h^{-1} : V^{2} \rightarrow u=v$ analogous with (\cite{BrelotChoquet51}), where $h^{-1}$ preserves regularity, we have that
$f \in C^{\infty}(\Omega \backslash V^{2})$ iff $(f- \delta) \circ h^{-1} \in C^{\infty}(\Omega)$, for some  $V^{2}$.
T has point support, when it is extremal (\cite{Choquet62}), that is T(u)=T(v) implies u/v=c. Assume for instance $f-\delta \in C^{\infty}$, with supp $\delta \subset$ the representation space and $supp \delta =h^{-1}(V^{2})$, then we assume supp $\delta$ regular, according to (\cite{BrelotChoquet51}).

\subsection{Contraction domains}
\label{sec:contraction}
Assume $f=e^{g}$ and $g=e^{\phi}$, with $\mu(f)=\widehat{\mu}(g)$, for a $\mu \in \mathcal{E}^{'(0)}$.
Relative an involution condition, it is sufficient to consider the phase space and $f_{j}(u)(x) \rightarrow f(x)$ defines ``weak'' convergence. $\phi=const \neq 0$ implies $f=const$, assumes an involution condition. Example: df=dg f= $d \phi g f$. Given $\phi$ absolute continuous, we have $d \phi=0$ iff $\phi=const$ implies df=0. Further, where $\phi=const$ implies $f=const$, f is absolute continuous.
Assume $\mathcal{L}^{2}I(\phi)(\xi + i \eta) \rightarrow \mathcal{L}^{2}I(\phi_{0})(\xi + i \eta)$, with $\phi \in H$, $\xi \in \Gamma$ convex and majorized by polynomials in $\eta$. Then the convergence for $g \in C^{\infty}(\mbox{int} \Gamma)$ (inner points) follows. By completing g to H and majorized by polynomials, repeated completion induces convergence for f in $C^{\infty}(\mbox{int } \Gamma)$. Completion corresponds to approximative solutions. Note that the completion affects the order of the group for approximative solutions.

\vsp

Example: dF(x,y)=$F_{x} d x + F_{y} d y$ a $C^{\infty}(dx,dy)$ form, does not imply a $C^{\infty}(dx)$ form. dF=0 implies
$F_{x}/F_{y}=-dy/dx$, that defines a contraction domain, when $F_{x}/F_{y} \neq const$, that is monogenity
is dependent of $x \rightarrow y(x)$ and $dx/dy \rightarrow dy/dx$.
Contraction according to $dy=\sigma(x) d x$, with $\sigma \in C^{\infty}$, defines a $C^{\infty}$(dx) form.
 Given dF=0, that is a closed form in (dx,dy), under the condition $F_{x}/F_{y} \neq const$, we can determine y, such that dF is a closed form in dx. On a contraction domain, $\sigma= dy/dx =-F_{x}/F_{y}$ and when the regularity is preserved under $\sigma \rightarrow 1/\sigma$ , closed forms are preserved.

 \vsp

Consider $G_{s}$ according to $U_{s}f=c$, when f=c, that is movements that preserve constant surfaces. Assume $G_{s}$ has range in the polar set.
A regular approximation property implies abscense of invariants. When Uf=c implies U=I, we see that the approximation property is relative type.

In analogy with the implicit function theorem,
assume $\Phi  : \quad G \times G \rightarrow G$ a $C^{1}-$ mapping, with $d \Phi/dU$ invertible, then we can determine a sub-group $G' \subset G$, such that for every $dV \in G'$, we have unique existence of $dU \in G$, with $\Phi(U,V)=0$ and $\frac{d \Phi}{dV} \frac{d U}{d \Phi} dV=dU$.
The contraction domain is assumed dependent on the boundary, that is $\Phi=0$ at the boundary. When $\Phi$ is only continuous , where G is defined in $\mathcal{E}^{'(0)}$, but defines continuation according to a normal model, then dU can be determined as a regular (weak) limit, when $dN \rightarrow dI$.

\vsp

 Assume f continuous and  real. For the transformation $Uf=\int f d U + \int_{\alpha} df$,  we start with $d U \rightarrow (dU,dV) \rightarrow dN$ and consider $\int f (dU - d U_{N}) =\int_{\alpha} df$,  where the boundary is given by dN/dV=0, that is $d U_{N}=dN(U,V)$, where dN a normal operator (\cite{AhlforsSario60}).
 Example: when $dI \in G$ and $f \in \mathcal{H}(G)$, with regular approximation property, we can define a normal model.
 A necessary condition for existence of a normal approximation N at $\alpha$, is that Flux(Nf)=0 over $\alpha$. Assume $\int_{\alpha} df^{\diamondsuit}=\int_{\alpha^{\diamondsuit}} df$. When f is analytic, we consider $f=e^{\phi}$, with $\phi$ completed to harmonicity, $d \tilde{U}_{N}(f)=dU_{N}(\tilde{f})$ and $\tilde{\alpha}$ thus takes into account $\tilde{\phi} - \phi$. That is we assume $\int_{\alpha} d \tilde{U}_{N}(f)= \int_{\tilde{\alpha}} d U_{N}(f)=0$. When we consider completion to $L^{1}(G)$, it is sufficient to consider
 $\mid f \mid (u,v) \in L^{1}$.

Starting from $\int f d U + \int_{\Gamma} df$, we see that $Uf-f=\int f (d U - d I) + (\int_{\Gamma} - \int_{0}) d f$. Monogenity implies $\lim_{\Gamma \rightarrow 0} f'(x)=f'(0)$, that is sufficient  for $\Sigma=\{ dU=dI \}=\{  U=I \}=\Delta$, is $\int_{\Gamma \backslash 0} df=0$.
Example: $\int_{\alpha} df=f(0)$, that is $\alpha \sim 0$,  where $\alpha$ defines a single valued continuation. Absolute continuous topology, is interpreted as the continuation being independent of the boundary.

  \subsection{Normal surfaces}
\label{esec:normal}
   Single valuedness for the normal assumes  that $d U \times dV$ generates surfaces, that is we assume $d U \not\equiv dV$. Multivalence relative conjugation $d U \rightarrow dV$,
   is considered as multivalent continuation and monogenity gives locally unique continuation. We assume
 $f_{v}=\rho f_{u}$, for $\rho$ regular, defines a contraction domain. That is f is dependent of both u,v outside the boundary, corresponding to a regular approximation property. Volume preserving conjugation is according to $\rho d I \rightarrow (1/\rho) dI$ absolute continuous and defines a contraction domain, when $\rho$ regular and finite.

  \vsp

  Assume $Uf=\mathcal{L}^{2}I(U \phi)$, with $f=e^{g}$ and $g=e^{\phi}$ analytic. Given $<g,dV>$ linear for $d V=\rho dI$ and $\rho g \in L^{1}$, we have that $\tilde{V} e^{g}=e^{Vg}$ has algebraic singularities. Given Vg analytic, bounded and  single valued in the plane, it is linear. Consider for this reason $f=\mathcal{L}I(g)(u,v)$, for instance $Uf=\mathcal{L}I(Vg)$. Then under the condition that $<g,dV>$ is linear, an approximation property for dU, is independent of dV.

  A boundary in $G_{3}$ is 2-dimensional, that is $G_{2} \subset G_{3}$. G continuous implies $dI    \sim \Sigma dU_{j}$, with $d U_{j} \in G$. In the representation space, $G_{2}$ continuous, defines a plane and thus determines two convex  sets in $G_{3}$. Given $f \in H(G_{2}=bd G_{3})$ implies $f \in H(G_{3})$, we have a complete max-principle.
  Traces are regarded as irreducibles in representation space. Example: volume preserving conjugation : $\rho \rightarrow 1/\rho$ absolute continuous, with $(dU,dV)=(1/\rho,-\rho)dI=-(\rho,1/\rho)^{\diamondsuit}dI$. Note that $(\rho + 1/\rho)^{-1} \rightarrow 0$, when $\rho \rightarrow 0$ or $1/\rho \rightarrow 0$, that is we have a strict condition for $d U \times dV$, also when we have a non-strict condition for dV.

\vsp
Assume $\phi : f(x) \rightarrow f(1/x)$ algebraic, then constants are preserved. $-d f(x_{0})=x^{2}d f(x)=0$, for $x_{0}=1/x$.
For the condition $x_{0} \rightarrow 0$ to preserve a rectifiable mapping, we must have density for translates. Example: $1/ \mid x-y \mid \rightarrow 0$ also when $x \neq y \rightarrow \infty$, that is we have singularities outside traces.
Example: $f(x)=f(y)$ according to $f(x/y)=f(1)$ implies
$df(x/y)=0$ implies $x/y=const$, that is point support for dI implies an absolute continuous mapping $(x,y) \rightarrow x/y$.

$\overline{N(f(x) f(1/x))}=\overline{N(f(x))} \cap \overline{N(f(1/x))}$, simultaneously $N(f(x) f(1/x))=N(f(x)) \cup N(f(1/x))$, given that $1/f(x) \neq f(1/x)$, for instance $f (x) \rightarrow f(1/x)$ algebraic.
The first equality implies irreducible first surfaces relative $x \rightarrow 1/x$ or restriction to decomposable first surfaces relative $x \rightarrow 1/x$. Example: x,1/x are regularly situated.

\subsection{A strict condition}
\label{sec:strict}

Assume $\eta_{1},\eta_{2}$ give two classes of transversals, according to
$<\eta_{1},d {}^{t} U>(\phi)=<\eta_{2},d {}^{t} V>(\phi)$. Given $d {}^{t} V=\rho d {}^{t} U$, we see that
$<\eta_{1},d {}^{t} U>(\phi)=< \rho_{1} \eta_{2}, d {}^{t} U>(\phi)$, where $\rho_{1}$ the transponate relative $<,>$.
Assume for this reason $\eta_{1} \rightarrow \eta_{2}$ does not preserve type. Example: $\eta_{2}(x)=\eta_{1}(\phi(x))$,  where $\eta_{1}$ algebraic (strict condition) and $\phi$ absolute continuous . Change of type is necessary for the polar set to be given by
a 2-dimensional surface, $\eta_{1},\eta_{2}$. This however does not imply regularity.

\vsp

Example: consider $\eta_{2}(x)=\eta_{1}(\phi(x))$, with $\eta_{2}/\eta_{1}=1/q(x)$ and $q(x)$ polynomial. Then, when q is reduced, we have $\eta_{2}/\eta_{1} \rightarrow 0$, that is $\eta_{2} \bot \eta_{1}$ (pseudo-orthogonal). Further, we have $q'/q \rightarrow 0$ regularly. $\phi'(x)=1$ does not correspond to a strict condition. When $\phi(x)$ is assumed to separate convex  sets in the plane, this implies presence of irreducible (linear) boundary.

Assume $h \in H$ and $f,g \in L^{1}$. Uf=g and Vh=g, implies $dV^{-1}U \in G_{H}$, that is $\rho f, \vartheta h \in L^{1}$ and $f=(\vartheta / \rho) h$.
In particular, when $\vartheta=1/\rho$, we assume $\rho f \in L^{1}$ iff $h/\rho \in L^{1}$, for instance when $\rho=$ constant.
When $d \rho \rightarrow 0$, we have $d (1/\rho) \rightarrow 0$, as long as $\rho^{2} \rightarrow const$. Sufficient is that dU is algebraic over $L^{1}$. Any $G_{H}$, has a 2-parameter subgroup. Consider $\varphi_{1} d U_{1} + \varphi_{2} d U_{2}=\omega d I$. Assume $d W_{j}=(\varphi_{j}/\omega) d I$,
then $d W_{1} + d W_{2}=dI$, that is a projective decomposition.
We have  according to Lie four cases (\cite{Lie91}): $\omega \equiv 0$ or $\omega \not\equiv 0$ and $\big[ dW_{1},dW_{2} \big] \equiv 0$ or for instance $dW_{1}$.
Necessary for a projective decomposition, is that $\varphi_{j}/\omega \in L^{1}$, for instance $\varphi_{j} \prec \prec \omega$.
Example: assume $\rho d U_{1} + (1/\rho) d U_{2} = \omega d I$ with $\rho/\omega,1/\rho \omega \in L^{1}$,
when $\omega=\rho + 1/\rho$, $\rho/\omega \rightarrow 0$, as $\rho \rightarrow 0$ and $\rightarrow 1$, as $\rho \rightarrow \infty$. Simultaneously $(1/\rho)/\omega \rightarrow 1$, as $\rho \rightarrow 0$ and $\rightarrow 0$, as $\rho \rightarrow \infty$. Necessary for $1 \in L^{1}$ is that  the support is bounded. Assume $\rho \rightarrow 1/\rho$ absolute continuous and $dI \in G$, defines extremal rays to convex cone. Note that a cone domain, is not a nuclear space. Thus, it is necessary that $\rho \neq const$,
that is we have change of type of movement.
Note that the order of a convex surface defined by G, can be defined by the number of invariants that limit the surface. Example: a conical spiral is limited by three invariants, a cylindrical spiral by two.
Note also that for multivalentness with convex  planar leaves, when the normal combines the leaves, an orientation is necessary, to separate between inner and outer domains.

\vsp

Assume for  real x,y, $B(x,y)=\Sigma_{i} c_{ij} x_{i} y_{i}$, a bilinear form, for positive scalars $c_{i}$. This implies that $B(x,x)=0$ implies $x_{i}=0$.
Given $B(x,y)=1$, $x \rightarrow y$ defines duality with respect to the scalar product. In particular for a line L,
$t x_{i} \in L$ iff $(1/t) y_{i} \in L$, for $t >0$ real.
Example: $\alpha + \beta=\alpha \beta$, given $\alpha \beta =1$ implies $\alpha=1$ or $\beta =1$, thus
 we must have $\beta=0$ or $\alpha=0$, thus $\alpha \beta=0$!, that is a projective dual conjugation is reducible.
 When (x,y) gives a domain for absolute continuity for B, then
$d B=0$ with $\Sigma c_{j}y_{j}d x_{j}=0$ corresponds to $d I \in G$, where $K(G)=G \times G$ is defined  with respect to B=1.
 \vsp

Given translation is dense in (u,v), we have (u,0) + (0,v)=(u,v).
Example: $f(u,0)=f_{1}(u)$ and $f(0,v)=f_{2}(v)$·
u+v is dense on compact sets in the plane, given $f(u,v) \in \dot{B}$. Assume $\mid f(u,v) - g(1/u,1/v) \mid < \epsilon$,
 where $g \rightarrow 0$, when $u \rightarrow 0$ and v fixed, in the same manner, for u fixed and $v \rightarrow 0$. That is, $g \in \dot{B}(nbhd \infty)$ and $f \in \dot{B}(nbhd 0)$, when we have that the translates are dense close to 0.
Assume for $\mid \alpha \mid > k_{0}$, $D^{\alpha}f=0$ implies x=0, that is a strict condition. Define boundaries $\Gamma_{k}=\{  D^{\alpha}f=0 \quad \mid \alpha \mid \leq k \}$, these define a very regular boundary. $F_{\alpha}(x,y)$ linear implies that sub-level surfaces contain disks, that is of order 0. $UVF_{\alpha}=VU F_{\alpha}$ regular, when $V \rightarrow I$, corresponds to algebraicity for U. $dI \in G$, with $d V \rightarrow dI$, over a contraction domain, gives reduction to an algebraic boundary. The approximation property in U, gives algebraicity in (U,V) as above (cf. section \ref{sec:nuclear}).
Example: $R(dU)^{\bot}$ closed, does not imply that $H(R(dU)^{\bot})$ is nuclear.
Assume $R(dU)^{\bot}$ not nuclear, that is does not have the approximation property. Note that existence of dW, that gives completion through a polar set to nuclearity, does not imply regularity.
Consider continuation through dV, such that the diagonal u=v is a closed set. When $R(dU,dV)^{\bot}$ is given by R(dN), the condition that dN is independent of V at the boundary, is necessary for a regular normal approximation property.

\vsp

In order to discuss T extremal on a space X, we assume the inclusion is related to $(I_{1}) \subset (I_{2})$, that is a majorization condition that implies a strict condition.
Assume $T(f) \leq m$ on leaves $L_{j} \subset X$ and a complete max-principle on X. Then T(f) is bounded on trajectories between $L_{j}$.
Consider $T : (I_{2}) \rightarrow (I_{1})^{\bot}$ as restriction. Given $T : (I_{2}) \rightarrow (I_{2}) \backslash 0$ is surjective, we have a trivial restriction. Example: $(I_{2})=R(dV)$
and $(I_{1})=R(dU)^{\bot}$, with (dU + d W) + dV=dI and  where $dW \neq 0$, we may have a diffuse front, when $dW \rightarrow 0$ (cf. section \ref{sec:monogenic}, partial monogenity).
Starting from (\cite{Riesz27}]) and $\int f(u,v)du dv=\int f(x) d \varphi_{1} \otimes d \varphi_{2}$,  where $d \varphi_{1}$ and $d \varphi_{2}$ have bounded variation. Assume $\Omega_{j}=$ $\{ (u,v)_{j} \quad f((u,v)_{j})=const \}$, that is first surfaces corresponding to a 2-parameter subgroup of G. Further, $\Phi : \Omega_{j} \rightarrow \Omega_{j+1}$ absolute continuous, $j=1,2,\ldots$. Then we have that $\Phi$ corresponds to a zero-function. We assume $\{ (u,v) \}$ dense in the domain for f, such that  $f \in C^{0}(x)$
can be approximated by $H((u,v))$. Zero-sets are defined by V such that  $\mid d (\varphi_{1} \otimes \varphi_{2}) \mid < \epsilon$ (total growth). In particular, given volume preserving measures according to $\rho d I \times (1/\rho) dI=dI$ on a planar domain, we assume that zero sets have measure zero.
$\Omega_{1} \rightarrow \Omega_{2}$ has a  corresponding mapping $d (\varphi_{1} \otimes \varphi_{2})_{1} \rightarrow d ( \varphi_{1} \otimes \varphi_{2})_{2}$ absolute continuous and thus maps sets of measure zero on sets of measure zero.
Consider $f \in L_{\varphi}^{a}$ according to $\int \mid f \mid^{a} d \varphi$ finite. For a linear functional B, $\mid B(f) \mid \leq M \parallel f \parallel_{a,\varphi}$, we have existence of a generatrix $w \in L_{\varphi}^{a'}$,
except for a zero function,
with $M=\parallel w \parallel_{a',\varphi}$ and $1/a+ 1/a'=1$ (\cite{Riesz27}).
Example: assume $d (\varphi_{1} \otimes \varphi_{2})_{2}=\rho d (\varphi_{1} \otimes \varphi_{2})_{1}$, for $\rho$ regular in $L^{1}$, then $<f(w),d (\varphi_{1} \otimes \varphi_{2})_{2}>=<{}^{t} \rho f(w), d (\varphi_{1} \otimes \varphi_{2})_{1}>$. Let $\rho_{1}$ be the transpose, relative $<f,w>$, then we have:
\newtheorem{lemma55}[lemma1]{Lemma}
\begin{lemma55}
Given $T : L_{(\varphi_{1} \otimes \varphi_{2})_{1}}^{a} \rightarrow L_{(\varphi_{1} \otimes \varphi_{2})_{2}}^{a}$ continuous, there is a continuous mapping between the corresponding generatrices.
\end{lemma55}

\subsection{A strict subgroup}
\label{sec:subgroup}
Consider a strict subgroup according to $d U \bot d P$, where d P algebraic (normal).
Assume sng is defined by traces (irreducibles), that is we assume $f(x,y) \in C^{\infty}$ $x \neq y$ and f very regular. Given density for translates, the polar set is trivial. Otherwise, consider regular approximations  $(x+w,y)$, that is where regularity is dependent of w. A strict condition implies independence  of the parameter w.

$dN \in C^{\infty}_{c}$ implies locally existence of $d N_{j} \rightarrow dN$, with $d N_{j}$ algebraic. Example: The proposition f(u,v,n) regular implies f(u,v) regular, corresponds to regularity being preserved, when $dN \rightarrow dI$.
Consider $d U \rightarrow dV$ as continuation, $d (U,V)^{*}=-d (U,V)^{\diamondsuit}$, with $V=U^{*}$. When given $d (U,V) \rightarrow d (U,V)^{\diamondsuit}$ projective, then $d (U,V)$ is harmonic and $d(U,V) \rightarrow d (U,V)^{*}$ preserves regularity in the plane (\cite{AhlforsSario60}).

Assume f hypoelliptic with algebraic polar, then $f(G_{r})$ reduced, can be continued to $f(G)$ reduced. Necessary for monogenic continuation through v, is that v=v(u) and $d V/d U \rightarrow 0$ regularly at the boundary. It is further necessary, that the domain can be reduced to the boundary. The boundary to a pseudo convex domain is of order 0 and when (u,v) generates a plane, it has 1-dimensional normal in 3-space. Conversely, a 1-dimensional boundary in 3-space, implies a non-regular (multivalent) normal.

Consider $(U,V)f \rightarrow 0$, when $f \rightarrow 0$. The order of zero's can be determined by determining the order for restrictions to complex lines, why it is  sufficient to consider pluri complex representations.
Assume $G(dU)=\{ d V \quad UVf=VUf \rightarrow 0 \quad f \in (I)_{B} \}$, then pseudo topology according to $\dot{B}$ (\cite{Schwartz66}), gives a representation in $\mathcal{D}_{L^{1}}'$,  where $(\dot{B}')'=\mathcal{D}_{L^{1}}$, that is we assume for $d U \in G$, when $\phi \in \mathcal{D}_{L^{1}}$, that $<T,dU> \in \mathcal{D}_{L^{1}}'$, $ \forall T \in \mathcal{D}_{L^{1}}'$. Example: $\mathcal{L}(\phi)=gf-1=\mathcal{L}(G*F) - \mathcal{L}(0)$, that is outside the diagonal $\phi=G*F$ in $\mathcal{S}'$, the conjugation is taken according to $UG*F=G*VF$, convergent outside the diagonal.

\vsp

Example: starting from a neighborhood, defined  by $d \mid F \mid < C$, with  $d>0$ minimal (geodetic), we have that $\mid F \mid < C/d$, that is we have a complete max-principle. Given analyticity, max is not reached in the neighborhood.
In this case, presence of a max-principle, implies presence of a polar set.
Assume $d \mid \phi \mid < C$, for $F=e^{\phi}$ and when $d=0$ corresponds to the polar set to F, outside the polar set, F can be given a regular representation.

Consider a normal model, $dN \bot dU \times dV$, assume $dN \bot dV \times dU$. Necessary for $dI \in G_{H}$, is that  dN is uniquely determined by $dU,dV$. Sufficient is that $f_{u},f_{v}$ regular are linearly  independent. dN is defined as in (\cite{AhlforsSario60}) for $f \in C^{0}$
and flux(Nf)=0, that is $\int_{\Gamma} d N(f)^{\diamondsuit}=0$. Assume $(dx,dy)^{\diamondsuit}=(-dy,dx)$ and $d f=\alpha dx + \beta d y$,
then we have $-d f((dx,dy)^{\diamondsuit})=df^{\diamondsuit}(dx,dy)$, that is $\int_{\Gamma} df^{\diamondsuit}=-\int_{\Gamma^{\diamondsuit}} df$.

Example: consider $M(\vartheta,\vartheta^{-1})$, with $\vartheta^{-1}=x/y$. We have that $d \vartheta^{-1}/d \vartheta=-1/\vartheta^{2}$ is single valued as a function of $\vartheta$, when $\vartheta$ is real. In particular,
$\frac{d M}{d \vartheta^{-1}}=\vartheta^{2} \frac{d M}{d \vartheta}$ and given M monogenic in $\vartheta$, then M is monogenic in $\vartheta^{-1}$.

\newtheorem{prop1}[lemma1]{Proposition}
\begin{prop1}
 A normal model induces a monogenic subgroup.
\end{prop1}

Assume existence of dN(U,V), with $dN/dV \rightarrow 0$ regularly at the boundary.
The condition dU=dI implies dV=0, with $dV \in G$, defines regular invariants.
For a monogenic subgroup of $G \times G$ defined by conjugation, according to (u,v), we have that  $f'(u,v)$ is continuous,
implies $f'(u)$ continuous, that is a regular approximation property in (u,v) implies a regular approximation property in u.
For $Uf=\int f dU + \int_{\Gamma} df$, we have Uf=f iff $\int f d (U -  I) +\int_{\Gamma} df -f(0)=0$.

Example: assume f(u,v) regular and $\lim_{v \rightarrow 0}f(u,v)$ regular, according to $\lim_{v \rightarrow 0} \frac{\delta f}{\delta u}(u,v) \neq 0$.
When $\vartheta_{j} = \delta u / \delta x_{j}$, such that  $\vartheta_{j} \frac{\delta f}{\delta u} (u,v)=(U,V) \frac{\delta}{\delta x_{j}} f(x)$ continuous, a projective (orthogonal) conjugation  (U,V)=I is sufficient for a normal model.

Example: assume $d U^{\bot \bot} \rightarrow d U^{\bot}$ not projective, but for instance $d U^{\bot \bot} / d U^{\bot} \rightarrow 0$ in $\infty$ regularly.
When $d U^{\bot}$ is completed to $L^{1}$ (\cite{AhlforsSario60}), that is $R(dU^{\bot})$ is closed, $d U^{\bot \bot}$ can be defined by annihilators.
Assume $R(dU^{\bot})$ is closed and $N(dU^{\bot \bot})={}^{\circ} R(dU^{\bot})^{\circ}$, thus $dU^{\bot \bot}(f)=0$ implies $f=dU^{\bot}(g)$, for g in the domain to $dU^{\bot}$.

Fix dU and define dV through UV=VU. R(UV) regular, defines a domain for convergence.  Assuming Uf of negative exponential type, when VUf changes exponential type, we have singularities.
Example: $\log \mid Uf \mid \rightarrow \alpha \mid x \mid$, with $\alpha<0$ and $\log \mid Vf \mid \rightarrow \beta \mid x \mid$, that is when $\alpha + \beta >0$, the movement changes exp-type.
exp-sgn(u,v)=exp-sgn (u) $<0$ defines a very regularizing subgroup.

\subsection{The involution condition}
\label{sec:involution}
Given $\alpha + \beta \in sng$ implies $\alpha,\beta \in sng$, we have that sng defines an extremal ray in a convex cone.
Example: $g=e^{C_{2} \phi}$ implies $g'(x)=C_{2}g$ and we have existence of $C_{1}$ with $C_{1} g'/g =const>0$, that is we can determine $C_{1}$, such that
$C_{1}C_{2}>0$. In particular, $\Sigma C_{1,j} g'=\Sigma C_{1,j} C_{2,j} g$ and $\Sigma =0$, implies $C_{1,j} C_{2,j}=0$ $\forall j$.  $g'=const$ represents sng relative monogenity, that is we do not have a contraction domain.

\vsp
 Assume P a convex cone in E and $E=F+P$, then a continuation to $O \supset P$ is unique, if $\sup_{a \leq x} \tilde{f}(a) = \inf_{x \leq b} \tilde{f}(b)$ (\cite{Choquet62}).
 Assume F+P corresponds to $R(dU)_{H}$, that is the range for analytic action of movement. Thus, $R(G)_{H}$ is represented by $F + \cup P_{j}$,  where the order for G determines the number of convex $P$. Since $R(G_{2})_{C}=R(G_{8})_{H}$ in the plane, the order is dependent of the topology.
Note that scaling of $x_{1},x_{2}$ gives $\xi_{1}/\xi_{2} = - (t d x_{2} / t d x_{1})$, that is every involution condition defines a convex cone.
Further, $d \xi/ d \eta= \frac{\xi_{x} dx + \xi_{y} d y}{\eta_{x} d x + \eta_{y} d y}$. Thus. we have $\frac{d \xi}{d \eta}(tx,ty)=\frac{d \xi}{d \eta}(x,y)$.

Assume for a measure $\mu_{t}$, $T=  \lim_{t \rightarrow 0} \frac{1}{t} (\mu_{t} - \delta)$,
with an approximations property of type (H) (\cite{Faraut70}), then T defines uniquely a measure $\mu_{t}$, that identifies type of movement.
Assume approximation through dU=0 according to
$\xi d x + \eta d y=0$, independent of t. $d U \rightarrow dI$ corresponds to $\xi,\eta \rightarrow 1$ under preservation of the involution condition. Example: $f \in L_{ac}^{1}(\Gamma)$, with dU(f)=0 implies U=I,
that is $\int_{\Gamma} df=f$ over $\Gamma$. When $\mu_{t}$ is very regular, we have $\mu_{t} : \mathcal{D}' \rightarrow \mathcal{D}^{' F}$, in particular there is a corresponding polar set in $\mathcal{D}'$,
that is when $\nu_{t} \mu_{t}=\delta$ in $\mathcal{D}'$, then ker $\nu_{t}$ is non-trivial.
Example:  When the domain for f, $\Omega$ very regular in the sense of (\cite{BrelotChoquet51}), this implies parabolic boundary (\cite{Parreau51}), that is $f \in R(G)$ has exponential representation.

\subsubsection{Resolution of sng}
We consider three types of singularities, $\Omega_{0}:$ according to $x \in \Omega$ not accessible for $dU=\Sigma dU_{j} \in G_{r}$, that is
$\Omega_{0}$ represents the polar set,
for instance $d I \notin G_{r}$.
$\Omega_{1}:$
according to $x \in \Omega$, but with abscense of disk subset, for instance abscense of a very regular subgroup of $G(\Omega)$.
Thus, $\Omega_{1}$ defines singularities relative $C^{\infty}$.
Example: $t^{2}f(x,y)(t)$ gives two -valent continuation, that is
a convex domain in two layers, gives one-sidedness relative two lines in the plane and we can have abscense of disk subset.

Consider also $\Omega_{2}$ according to
$x \in \Omega$ not accessible linearly, that is through one-parameter approximation.
Thus, $\Omega_{2}$ defines singularities relative $G_{H}$, that is relative a continuous group. Example: dN(U,V) is dependent of both arguments at the boundary, that f is not monogenic or $x \notin $ disk and $f(u,v) \notin C^{\infty}$
Finally, note that monodromy excludes traces, for instance $d V \rightarrow dI$ excludes $d U \rightarrow dI$. Consider reduction to H and to an algebraic model relative projectivity, that is $d(U,V)=d(U,I) \otimes d(I,V)$ corresponding to a decomposable topology. Note that (J) is a maximal domain for regularity iff (J) is irreducible. If $(J)=R(dU,dV) \subset X$, can be continued regularly through dW, given (J) maximal, we must have dW=dI. However, polar set can be given by non-regular continuations.

\vsp

 Assume X the domain for $f \in \mathcal{H}$ and consider a subgroup $G_{r} \subset G$, according to $G_{r} X \subset X$. We want to determine $d W \notin G_{r}$,
with $dW X \subset X$, that is dW in the polar set in X relative $G_{r}$. We assume further that the continuation is regular and such that  UW=WU, that is (dU,dW) defines a very regular subgroup of G,  where $d U \in G_{r}$. We also assume maximality in the following sense: $dW/d U = \rho$ and we do not have existence of $d V \in G_{r}$, with $dW/dV=\sigma_{1}$, $dV/dU=\sigma_{2}$ and $\rho=\sigma_{1} \sigma_{2}$ regularly. When the continuations are assumed regular, we exclude presence of traces.
The co-dimension here denotes the smallest number of $d W_{j}$, linearly independent, such that $d U + \Sigma d W_{j}=dI$ in H.
Thus, $d U \rightarrow dI$ implies $\Sigma d W_{j} \rightarrow 0$. Note that non-existence of dV between $dU \rightarrow dW$ and volume preserving, according to $\sigma_{2}=1/\sigma_{1}$, defines a maximal chain.
Example: Assume F a first surface, multiply connected, then the number of regular normals (linearly  independent) defines the co-dimension.

When $dI \notin G_{H}(j=2)$, but $dI \in G_{H}(j=8)$, we can determine a 2-subgroup in $G_{8}$.
Example: $f(u,v)=0$, for some (u,v) defines a path. Assume further $dI \in G_{C^{\infty}}(j=2)$,
but $d I \notin G_{H}(j=2)$. $d I \in G_{H}(j=8)$ implies $d I \in G_{C^{\infty}}(j=8)$ assuming a strict condition. We have  a contraction domain, given abscense of invariants. Example: $d u_{j}/d u_{1} = \rho_{j} \neq const$ and regular with $df / d u_{1} \neq 0$, that is reduction to a conformal representation.

Assume starting from f(u) analytic, there is continuation along v in $C^{\infty}$, with (u,v)=(v,u) that is a disk,
where $(u,v)$ are related by conjugation. Example: $f(u,v) \in C^{\infty}$ implies $f(u,0) \in C^{\infty}$ with $dV \in G$. Sufficient  for this is monogenity.
Example:  $f \in H(G_{2})$, has continuation $\tilde{f} \in C^{\infty}(G_{2}) \simeq H(G_{8})$, given the domain a contraction domain, that is when $G_{8} \rightarrow G_{2}$ regular, we have analytic continuation. Assuming strict convergence, given q is accessible through $H(G_{8})$, we have accessibility through $C^{\infty}(G_{2}$).

\subsection{Exponential representations}
\label{sec:exponential}
Assume $L_{x}$ is defined by f such that $<f,{}^{t}P(D) \psi>=0$ and consider $X$ testfunctions corresponding to f. Assume $f \in \mathcal{D}_{L^{1}}'$ and for $D_{j}=\frac{\delta}{\delta x_{j}}$ and P(D)  a (homogeneously hypoelliptic) d.o. with constant coefficients, then given P(D)f=0, it is sufficient to consider test functions with representation in Exp on a compact set (cf. Montel theorem). More precisely, $\mathcal{D}_{L^{1}}$ with a strict condition, can be approximated by H, that is
$P(D)f(\psi)=0$, where ${}^{t} P(D) \psi$ is dense in $\dot{B}$. This is obvious when $\psi=e^{<\xi,x>}$, since $\{ {}^{t}P(\xi)=0 \}$ is trivial. When $\psi=e^{\phi}$ where $\phi$ single-valued and holomorphic on a compact set, it is sufficient to consider $\phi \sim <\xi,x>$.
When $\widehat{f(\psi)}$ is holomorphic in $\xi+i\eta$ and is majorized by polynomials in $\eta$ on compact sets in $\xi$, it is the Laplace transform to a distribution in $\in S'$ (\cite{Schwartz52}). $P(D_{x}) \mathcal{L}I(<\xi,x>)=0$ iff $P(\xi)=0$, that is $L_{x}=X^{\circ}$. Given the  domain for f is closed, ${}^{\circ} L_{x}=X$.
Assume $UV \mathcal{L}I(\phi)=\mathcal{L}I((U+V) \phi)$, assume $I \in \mathcal{S}_{x}'(\Gamma)$ (\cite{Schwartz52}),  where $\Gamma$ is an open convex set. Assume $F=\mathcal{L}I(\phi)(\xi,\eta)$, $\xi + i \eta \in \Gamma + i \Xi^{n}$, where F is majorized of polynomials in $\eta$ on compact sets in $\xi$. When $I \in \mathcal{D}_{L^{1}}'$, we have $\widehat{I}(\phi)=P(\xi) f_{0}$,
where we assume $f_{0}=\widehat{g}$, with g very regular, in particular $f_{0}$ locally bounded on a nbhd of the diagonal.
Example: when log F $\in L^{1}$ (algebraic singularities), a strict condition implies density for log F in H, implies convergence for I in $\mathcal{S}'(\Gamma)$.

\vsp

Given density for $R(dU,dV)$ in the domain for regular solutions, we have that $dV \times dU(f)=0$ implies existence of algebraic annihilators. Given a differential operator Q with constant coefficients, we have existence of
$P(D) e^{<x,y>} T=e^{<x,y>} Q(D) T$, for some differential operator with constant coefficients P. In particular when $T=e^{<z,y>}$, we have $e^{-<x,y>} P(D) e^{<x+z,y>}=Q(z) e^{<z,y>}$. When $z \neq -\infty$, we have $Q(z)=0$, for the homogeneous equation.

Assume $UF(\gamma)=F({}^{t} U \gamma)$, when $UF \in H'$ we have that ${}^{t}U \gamma \in H$ and we have existence of $d \mu \in \mathcal{E}^{'(0)}$, with $\mu=UF$. Example: $\gamma \in L^{1}$ with ${}^{t}U \gamma \in H$. Example: sharp fronts, relative
$F(\gamma)=const$, with $(dU,dV) \in G \times G$ very regular, implies $(U,V)F(\gamma) \rightarrow {}^{t}(U,V) \gamma \in C^{\infty}$ continuous.  Example: $F(\gamma)(u,v) \rightarrow 0$ in $\mathcal{D}_{L^{1}}'$, when $u \rightarrow 0$ for fixed v and ${}^{t}(U,V) \gamma \rightarrow 0$ in $\dot{B}$.

\subsubsection{Approximative representations}
By a suitable choice of regularizing terms, we can restrict $\mathcal{D}^{' F}$ to $\mathcal{E}^{'(0)}$. Sufficient for $Uf=<f,dU> \rightarrow d U \in \mathcal{E}^{'(0)}$ continuous, is that flux(f)=0 and $f \in C^{\infty}$. $T \in H'$ implies $T \simeq \mu \in \mathcal{E}^{'(0)}$.
Example: $f \in C^{0}(G_{2}) \simeq H(G_{8})$ and when $\Phi f=\int f d U$, $\Phi$ can be represented in $H'(G_{8})$.
Example: $T \in \mathcal{D}_{L^{1}}'$ implies $T=P(\delta) F_{0}$, with $F_{0}$ very regular.
In particular, there is a Fredholm operator F, with $F(T) \sim F_{0}$ (=0 modulo regularizing action). Assume $F({}^{t} P(\delta) \phi)=0$,  where $F=0$ is finite dimensional. Assume $T^{N}$ hypoelliptic, that is  corresponding to ker $F_{N} \subset C^{\infty}$, that defines a ``rectifiable'' localization.

\vsp

Given a semi-group of type (H): $T=\lim_{t \rightarrow 0} \frac{1}{t}(\mu_{t} - \delta)$ (\cite{Faraut70}).
Assume S(t) = $\mu_{t} - \delta$. N(S) defines the type of $\mu_{t}$. $S(t) \in C^{\infty}$ defines approximative solutions. $\lim_{t \rightarrow 0} S(t)=0$ in $C^{\infty}$, give a regular approximation property.

Consider $d U \rightarrow d (U,V) \rightarrow dN(U,V)$, over regularly situated domains (\cite{Malgrange66}), this assumes  dU=dV trivial (algebraic).
Example: $R(dU)=\Omega_{1}$ and $R(dU,dV)=\Omega_{1} \cup \Omega_{2}$, where $\Omega_{j}$ regularly situated, thus given $f \in C^{\infty}_{c}(\Omega)$,  where $\Omega$ is an open set, we have $\Omega_{1} \cap \Omega_{2} \subset \Omega$ has measure zero for higher order derivatives, that is $D^{k}f=0$ on a set of measure zero, for some $k \geq 1$.

\vsp
We define almost analytic functions (\cite{Garding87}): AE(f)=g with $\overline{d}g \equiv 0$ in y=0. Thus, monogenic extensions are almost analytic. AE(fg)=AE(f)AE(g). Thus $g \rightarrow 1/g$ preserves AE. Example: $d(1/g)=-dg/g^{2}=-(1/g) d \log g $, with $\log g$ monogenic for $g \neq 0$ single valued, that is given log g monogenic on single valued $g \neq 0$, we have that $g \rightarrow 1/g$ preserves monogenity.

\subsection{The maximum principle}
\label{sec:maximum}
 Consider (dU,dV) as a continuation, given dU maximal we have that (dU,dV)=(dU,dI) or f(u,v)=f(u,0).
dU=0 implies dV=dI gives a regular trivial continuation. When $R(dU)_{H}$ has an approximations property, dV=dI
defines a removable set.

Assume $\omega \in L^{1}$, that is $\omega dI \in G$. Given a strict condition, we have that $\omega$ can be approximated by H. Assume $\omega=e^{\phi}$, with $\phi$ algebraic. Then $\omega$ does not have essential singularities in $\infty$, that is $\omega=0$ is the limit of algebraic first surfaces.
We can determine $\varphi_{j} \in H$, with $\varphi_{1}/\varphi_{2} \sim \omega \in H$. In particular, when (u,v) open, we can write $\omega(u,v)=e^{\phi_{1}(u) - \phi_{2}(v)}$, for $\varphi_{j}=e^{\phi_{j}}$.

Consider a normal model dN(U,V),  where $d V/d U \rightarrow 0$ regularly
in $\mathcal{E}^{'(0)}$. Example: (P/Q + R)f=0 with Q reduced and R of negative type, implies Pf of negative type, that is $Pf + R'f=0$ corresponds to an inhomogeneous equation. P hypoelliptic implies f of negative type. For general P we consider a topological neighborhood of $C^{\infty}$. For a normal space, we have density  for $\mathcal{D}$.

\vsp

Assume $\Omega$ an open set, with reg $\Omega$ (regular points) 2-dimensional and that $f \leq m$ on reg $\Omega$ implies $f \leq m$ on $\Omega$. Example: $f(x) \leq m$ on $\mid x \mid > R$, then we have $f(1/x)f(x) =e^{\phi(1/x) + \phi(x)}$.
Assume L linear, such that $L(1/x) + L(x)=L(1/x + x)$ and such that $\phi$ is majorized by L.
Then, the condition $\phi-L \in L^{1}$, with L linear represents $f \sim P/Q$, where P,Q are analytic.
Example: $\mid f \mid \leq c \mid P/Q \mid$ implies $R(1/\mid x \mid^{k}) \mid f \mid \leq c$, where
$R(1/\mid x \mid) P(\mid x \mid) \leq c'$ implies $\mid f(x) \mid \leq P(\mid x \mid^{k})$ (\cite{Dahn23}).
When f is reduced,
we have that $f(1/x) f(x) \leq m$ implies $f(1/x) \leq m'$ using Nullstellensatz.

\vsp
Assume $\Omega_{m}=\{ F < m \}$ and $\Gamma=bd \Omega_{m}=\{ F=m \}$. When extremals are given by dF=0 on the boundary, given an absolute continuous topology, we have $F=m$ on the boundary. Thus, given $\Omega_{m} \subset \subset \Omega$, we have that the boundary can be given by dF=0 iff we have an absolute continuous topology.

Assume $A_{T}$ a domain for point support measure (\cite{Choquet62}). Then $a,b \in A_{T}$ and $f(a)=f(b)$ implies a/b constant.
Presence of disk-nbhd of 0 is necessary for the definition of $A_{T}$. Given a strict condition, such that  translates are dense in the plane, point support distributions define traces in the plane. When F is a very regular distribution, then $F \in C^{\infty}(x \neq y)$. Given X Hausdorff uniformity, with $u \rightarrow v(u)$ continuous, then u=v is a closed set (\cite{Bourbaki89}), thus $u \neq v$ is open, thus $H(u \neq v)$ is  decomposable, that is $\delta(u,v)=\delta(u) \otimes \delta(v)$. Example: $f \in C^{\infty}$ and $f \rightarrow c=f(0)$ does not imply $x \rightarrow 0$, when $\delta$ does not have point support, that is $A_{T}$ is a domain for a strict condition.
Example:
presence of continuum in the support for dI and $dU \in \mathcal{E}^{'(0)}$, implies at best $dU^{-1} \in \mathcal{E}'$, that is a discontinuous group.

\vsp

Example: fg-1 $\sim 0$ implies $g-1 \sim 0$ or $f-1 \sim 0$, given 1 irreducible (modulo regularizing action).
f hypoelliptic implies 1 irreducible (modulo regularizing action).

Example: $f \in C^{\infty}$ with $fg-1 \sim 0$ (modulo $C^{\infty}$). $f(x) \rightarrow f(0)$ does not imply $x \rightarrow 0$. $\delta_{x}(f) \rightarrow \delta_{0}(f)$ does not imply $x \rightarrow 0$, but given 1 irreducible, that is when f is reduced, we have monodromy in the sense of $x \rightarrow 0$.

\subsubsection{Parabolic surface}
Assume $\omega d I$ a harmonic measure, with
$\omega=0$ on the ideal boundary $\Gamma$, this implies $\Gamma$ parabolic, $\omega >0$ implies the ideal boundary hyperbolic.
For a unique normal operator dN, we consider a 2-dimensional surface $V^{2}$, note that flux is multivalent over $V^{1}$.

Presence of a max-principle corresponds to a parabolic surface, conformally equivalent with the plane $u + i v$, that is we have  abscense of traces. Example: u=v normal to a parabolic surface. The support to a very regular distribution has symmetric regularity outside u=v.
A very regular surface implies parabolic boundary implies parabolic surface.´

\vsp

When F is hypoelliptic, we have that the domain for convergence for G and FG=I is bounded, thus we have abscense of monogenity for G at the boundary, that is monogenity for parametrices does not imply monogenity for fundamental solutions. Example: FE=I + R,  where R is monogenous. That is we have that negative type for R(u,v) implies negative type for R(u).

Given sharp fronts, we have $\phi(u+v) \rightarrow \phi(u-v)$ preserves monogenity analogous with harmonic conjugation . Consider $f \rightarrow 1/f$, with
$f=e^{\phi}$ and $\phi$ odd. When $\phi \rightarrow \phi^{\diamondsuit}$ pure, $(x,y) \rightarrow -(-y,x)$ preserves analyticity in the plane, that is a two-mirror model according to $z \rightarrow z^{\diamondsuit} \rightarrow -z^{\diamondsuit}$, gives a pure model.
On a parabolic surface, $\phi<0$ implies $\phi=const$, thus continuation of a parabolic surface does not imply monogenity. When monogenity relative two mirrors can be reduced to monogenity relative one mirror, we have monogenity analogous with a normal model.

Presence of a separating functional, assumes inner points, $ v \neq 0$, that is it is necessary for monogenity, that we have
existence of $f'(u,v)$ single valued and finite, for $v \neq 0$. If further $d v/d u$ does not change sign, v can be defined as  monotonous in u.

\subsection{A boundary mapping}
\label{sec:boundary}

Starting from  a convex representation by G, given the boundary is given by dU=dI, we assume reduction to the boundary using a approximation property. $d I \in G_{\mathcal{H}}$ corresponds to nuclearity,
that is it is necessary for reduction to invariants, that the space is nuclear.

An analytic set A is irreducible iff reg A is connected (\cite{Chirka89}). Let $A=\{ d U =0 \}$, assume $d U \rightarrow dI$,
 where f=const corresponds to sng, according to $\int f d I=\int f'(x)dx$, that is  where dU analytic, we have that R(dU) convex implies reg A connected. When Uf is analytic has decomposition in irreducibles, it is  sufficient  to consider a pluri complex representation.

Assume $dI \in G$ extremal, gives reduction to dU=dI.  The domain =G can then from convexity be reduced to invariants, that is a multivalent  boundary. We consider contraction $G \rightarrow G_{2}$,  where $G_{2}$ is given by a conjugated pair and invariants $\Gamma_{i}$ are assumed rectifiable , that is we assume $\Gamma_{i} \rightarrow \Gamma_{i+1}$ absolute continuous.

\newtheorem{lemma4}[lemma1]{Lemma}
\begin{lemma4}
 When sng are represented by a continuous group $G_{H}$, we have resolution of sng by irreducibles.
\end{lemma4}

  Given R(dU), $R(dU)^{\bot}$ convex sets, assume p in an irreducible set in $R(dU)^{\bot}$ and dN continuation of R(dU) to p. When nbhd p is given by an analytic set, irreducibles can be chosen as lines.
  Singularities can be reached by irreducibles given co-dimension 1.
 Convexity is with respect to G, that is R(dU) are points accessible from G and $R_{H}(dU)$ gives regular points in R(dU).  A continuous group, has a decomposition in 1-parameter subgroups. Example: consider $G_{q}$ quasi inverses, that is $d U \rightarrow dU^{-1}$ does not change type. Then we have that $d UV=dI$ are irreducibles.

\vsp

Assume f=P/Q, for polynomials P,Q with Q reduced and consider $\frac{\frac{d Q}{d y} \frac{d y}{d x}}{Q} \rightarrow 0$ in $\infty$. Given sing supp Q $\varphi$ = sing supp $\varphi$, this implies that sng is determined by P(x), that is we have  monogenity. Over a contraction domain, such that $df/dn \neq 0$ and finite, the continuation  can be given by a normal (monogenic) model .

A family $\{ \gamma \}$ of continuous and rectifiable curves, of finite total length, is a normal family
(\cite{Oka60}) . A normal planar surface has a determined normal. Assume (u,v) defines a surface in the plane and
 $(d s/dt)^{2}=(du/dt)^{2} + (dv/dt)^{2}$,  where sng is defined by du=dv=0, in particular $dv \not\equiv 0$. When (u,v) defines a line v=cu, we have $(du/ds)^{2}=const$.  Single valued continuation in the representation plane according to $f_{u},f_{v} \neq 0$ linearly  independent, is mapped onto a single valued normal in 3-space. Example: continuation according to $u \rightarrow (u,v) \rightarrow n$, given regular contraction of formulas $f_{u}/f_{v}=n_{u}/n_{v}$.

 \vsp

The representation $Uf=\int f dU + \int_{\Gamma} df$, assumes that R(dU) $\subset R(U)$,  where $R(dU)^{\bot}$ is dependent of the boundary.
Necessary for a very regular approximation of surfaces, is that $\Gamma \rightarrow V^{2}$ a homeomorphism. Example: $(H_{-},H_{+})$ can also be seen as a two-valued covering of a disk (u,v), that is convexity with respect to two axes. Assuming codimension 1, $H_{-} \rightarrow H_{+}$ is defined by reflection through the boundary.

\vsp

Consider $\{ U_{j} f \}$, with $U_{j}f=\int_{(\Gamma_{j})} f d U_{j} + \int_{\Gamma^{\diamondsuit}_{j}} df$. Assume v=v(u), with $\sigma_{j}=dU_{j}/dU_{1}$ and $d U_{1} \sim d I$.  Then we have for sequential movements, $Uf =$ $ \int_{(\Gamma)} f dU + \int_{\Gamma^{\diamondsuit}} df=$ $(\int_{(\tilde{\Gamma}_{1})} \times \ldots \times \int_{(\tilde{\Gamma}_{N})} f du_{1} \ldots d u_{N}+$ $\Sigma \int_{\tilde{\Gamma}_{1}^{\diamondsuit}} \frac{\delta f \tilde{\sigma}_{j}}{\delta u_{1}} d u_{1}$.  Uf=f corresponds to multivalentness, for instance $d U=\rho dI$ with $\rho=const$ in $\infty$.
Example: $(d U + d U^{\bot})^{\bot}=\{ d U - d U^{\bot} \}$, that is $\sim 0$ on traces. Note that closedness with respect to $(d U,d U^{\bot})$ does not imply closedness with respect to $d U$.
Example: assume $d I^{\bot}= dV \in G_{H}$ defines an analytic set. Then (dU,dV) defines a Borelset, when $R_{H}(dU)^{\bot}=R_{H}(dV)$ (Suslin).
\vsp

Assume M closed, $M=\{ x \quad u(x)=0 \quad u \in L_{x} \}$. From annihilator theory $L_{x}=M^{\circ}=M^{\bot}$, when L linear. $L_{x}$ defines a general ideal (I)(M),  where $u^{N} \in (I)(M)$ defines rad (I). Assume $d U =\rho dI \subset T_{x}(M)$ (tangent space), for $\rho \in L^{1}$. Then, $L_{x}=\{ \rho(x)=0 \}$, that is given a strict condition, $L_{x}$ can be related to analytical ideals.

\subsubsection{Flux condition}

UV=VU implies $\int f (d UV-d VU) + (\int_{\alpha} - \int_{\alpha'}) d f$, that is given $\alpha'=\phi(\alpha)$, with $\phi$ absolute continuous, it is  sufficient  for UV=VU, that $\int_{\alpha} df-df \circ \phi=0$ , that is that $\phi$ preserves flux.

Assume df projective under $\alpha \rightarrow \alpha'$, that is $\int_{\alpha}=0$ implies $\int_{\alpha'}=0$.
Further, given decomposability, for instance $H(\Omega)$, when $\Omega$ open and contains $\alpha,\alpha'$, we have $\int_{\alpha}-\int_{\alpha'}=0$, implies $\int_{\alpha \backslash \alpha'}=0$.
Assume $Uf=\int f d U + \int_{\alpha} df=U_{0} f + \gamma_{\alpha}(f)$.
Then we have that
$VUf = \int f d (VU)_{0} + \int_{\alpha,\beta} df$. $VUf=\int f dU_{0} \times d V_{0} + V_{0} \gamma_{\alpha}(f) + \gamma_{\beta}(U_{0}(f)) + \gamma_{\beta}(\gamma_{\alpha}(f))$,  where we assume $V_{0} \gamma_{\alpha}$, independent of $\beta$ and $U_{0}$ does not depend of $\alpha$.
Example: sharp fronts implies $\int_{\gamma_{\alpha,\beta}}=\int_{\gamma_{\beta,\alpha}}$ and a unique continuation in $C^{\infty}$.

\vsp

When $\Gamma$ is multi-valent , we have multi-valent flux. Consider $\tilde{\Gamma}=(\Gamma_{0},\Gamma_{1})$,  where $\Gamma_{0} \rightarrow \Gamma_{1}$ continuous and  where $\Gamma_{1}$ satisfies a strict condition. In particular,
when $\pi: \Gamma_{1} \rightarrow \Gamma_{0}$ is absolute continuous, that is $d \pi(\gamma_{1})=0$ implies $\pi(\gamma_{1})=const$, then Flu$x_{\tilde{\Gamma}}$ can be defined  as single valued.
Example: assume G a non-parabolic surface and also $\pi^{-1}$ absolute continous.
Given $\Gamma_{1}$ 0-dimensional in G and $\Gamma_{0}$ 1-dimensional, using (\cite{BrelotChoquet51}) there is  a homeomorphism h, that composed with $\pi^{-1}$ non-trivial, maps $\Gamma_{0}$ on a continuation from the boundary, to inner points, that is starting from  points on web and an absolute continuous ``transversal'' continuation, the limit can be determined in inner points.
Example: scaling of a surface of revolution that defines a normal tube, that is the transversal intersects all boundaries and all in the same manner.

\vsp

Assume $\Gamma=\{ \Gamma_{i} \}$, i=1,2. a condition for monogenity on $\Gamma$, is then that $\lim_{\Gamma_{1} \rightarrow \Gamma_{2}} \int_{\gamma} d w= \lim_{\Gamma_{2} \rightarrow \Gamma_{1}}
\int_{\gamma} d w$.
$dw$ is closed if $\int_{\gamma} dw=0$ $\forall \gamma \sim 0$ and exact if $\int_{\gamma} dw=0$ $\forall \gamma$,  where $\gamma$ refers to closed contours. On a contractible space, we have $dI \sim 0$ (homotopical with the constant mapping). Example: assume $\Gamma_{1}$
with strict condition and $\Gamma_{2} \rightarrow \Gamma_{1}$ absolute continuous, corresponding to a contractible domain between the contours.
Assume $\int_{\Gamma} df=0$ and $\Gamma_{2} \sim 0$, for instance representing a trace that is not mapped onto sng f. Then the boundary mapping can be reversed and $df$ can be considered as exact. Otherwise
on the contractible surface that is limited by $\Gamma$, df is closed.

Concerning contraction of formulas, given $d \gamma_{1}/d \gamma_{2} \neq 0$ and f monogenous on a contractible domain, we have $\gamma_{1} \sim \gamma_{2}$ (homotopical).  Note that the normal dN, is not necessarily uniquely determined on $\Gamma_{2}$, but well defined on $\tilde{\Gamma_{2}}$.

\vsp

Example: $N(df^{2})=N(df) \cup N(f)$, that is a multi-valent boundary. Assume in particular the boundary is defined  by a strict condition, for some  iterate. Example: $P(D) \phi=0$ implies $\phi \in C^{\infty}$ defines a radical ideal $N(P) \subset C^{\infty}(X)$. Since monogenity is not a radical property, we can consider $f^{2} \in C(X)$ monogenous, for instance support decreasing by iteration, that is rad C(X) $\supseteq$ N(P).

\subsection{Projective conjugation}
\label{sec:projective}
Example: For $\sigma_{j} \in H$, let $dU + d V=\sigma_{1}dI + \sigma_{2} dI=dI \in G_{H}$,  where $G_{H}$ is continuous , with $\sigma_{1} + \sigma_{2}=1/p + 1/q=c$ iff $p+q=c pq$, that is we have a projective decomposition iff $dU+dV \equiv dUV=dVU$.  With these conditions, there is in H a maximal domain for projectivity. Further, $(dU,dV)^{\bot}$ defines the polar set that include traces.

\newtheorem{lemma6}[lemma1]{Lemma}
\begin{lemma6}
$G_{0}$ defines an extremal ray of a convex cone,  given
$dI \in G_{0}$ and $d U+dV=dI$ implies $dU,dV \in G_{0}$.
\end{lemma6}

Example: assume $d U,dI \in G$ and define
$G_{0}(dU)=\{ d U - d I \}$, then $G_{0}$ extremal means that $dI \in G_{0}$ implies $dU \in G_{0}$.
Assume $dU \in G$, then $G_{0}(dU)$ represents projective conjugates in G. When $dI \in G_{0}(dU)$, it is  a continuous subgroup of G. Further, when $d V \in G_{0}(dU)$, (dU,dV) forms a very regular subgroup iff U is algebraic, that is $dUI=dIU$.

\vsp

Assume $(U,V)f$ regular when $V \rightarrow I$.
dUI=dIU implies that $U \rightarrow {}^{t} U$ does not change type.
Assume $(dU + d V=dI)^{\bot}=\{ dU=dV \}$. Note that dU+dV=dI=dV + dU, that is projective conjugation  corresponds to a symmetric domain. However, we have that presence of traces implies presence of polar.
 Example: f(u,v)=0 with $d V/d U \rightarrow 0$, assumes a non-symmetric condition on the domain. However $(dU,dV)^{\bot}$ symmetric does not imply abscense of traces.

Example: assume hypo-continuity (\cite{Treves67}), for G according to $\pi(u-v) \leq p(u)q(v)$, with semi-norms in the right hand side. $\pi=0$ then defines a neighborhood of traces. Example: $v \neq 0$ with p(u)=0, further $\pi=0$, does not imply u=0, since p is a semi-norm.

\vsp

Assume UV is defined  (conjugation ) by $UV\mathcal{L}^{2}I(\phi)=\mathcal{L}^{2}I(UV \phi)$, sufficient  for this is that  U,V are algebraic (or UV algebraic). Given $g=e^{\phi} \in H$, with $\phi$ subharmonic. Consider $\tilde{\phi}$ completed to harmonic, such that $\tilde{g}$ is constant or analytic.
Further, $f=e^{g}$ with $<g,dV>$ linear and simultaneously $<\varphi,d {}^{t} V>$ linear, that is $dU \rightarrow dV \rightarrow d {}^{t}V$ linear ( single valued), conversely $dU=dU(dV,d {}^{t} V)$ and given canonical duality, dU can be defined as independent of $d {}^{t} V$.
Example: assume $\varphi \in X \subset \dot{B}$ can be determined such that  $X'$ can be given in Exp.
$\tilde{g} \in H$ and $\tilde{g}=\widehat{G}$, for some $G \in \mathcal{D}_{L^{1}}'$ (=$e^{-<\xi,x>} \tilde{\phi}$) implies that $\tilde{g}=\mathcal{L}I(\tilde{\phi})$ is possible.
Given X closed, we can approximate $X^{\circ}$ by Exp. Example: when $\varphi \in \dot{B}$ and $\widehat{g}=P(\xi) g_{0}$, where $g_{0}$ is of type 0, we can choose $\phi \in H$ so that $<g,\varphi> \sim <{}^{t} \mathcal{L}I(g),\phi>$. Modulo regularizing action, it is sufficient to consider P of type 0.

\vsp

$UF(\phi)$ absolute continuous is interpreted as F has a domain for absolute continuity in G. Example: $f \rightarrow 1/f$ continuous on an absolute continuous domain, that is local existence of fundamental solution. Volume preserving conjugation  can be compared to presence of fundamental solution.  Example: $UF*VG=F*G=I$ over (u,v), a domain for volume preserving conjugation. Regularity for $F(\phi)$ depends on the polar set. A strictly pseudo convex domain gives an algebraic polar set (transversal).

\subsection{Monodromy}
\label{sec:monodromy}
Borel (\cite{Borel12}) constructs a domain C=\{ x+iy \}, strictly larger than the domain for analyticity, for a monogenous function f.
C is related to a reduced domain $\Gamma \subset C$, its complement has measure zero in the plane, such that  when f is determined on $\Gamma$, f is determined on C. In particular, using the generalized Cauchy's formula, for $x$ on the boundary to $\Gamma$, f(x) can be determined as continuous from f(x+iy).
We will argue that for a monogenous function in the plane, $f \in C^{\infty}(x+iy)$ implies $f \in C^{\infty}(x)$,  independent of choice of y.

\newtheorem{prop2}[lemma1]{Proposition}
\begin{prop2}
Monogenity implies weak monodromy.
\end{prop2}

Consider $f(\phi) \sim 0$ (modulo regularizing action), for $\phi \in \mathcal{D}$ and f subharmonic,  with $f=\widehat{T}$ and complete $\phi$ to a harmonic function. Assume $L=\{ x \quad \phi \neq \tilde{\phi} \}$, then we can write $f(\phi)=<\widehat{T}(\tilde{\phi}), dI>+ \int_{L} d f(\phi)$. When f is monogenous, we can assume $\int_{L}=0$, that is it is sufficient to consider real $\phi$.
When T is of finite order, we can assume $\widehat{T}$ is majorized by polynomials, why it is sufficient
to prove $\mid \widehat{T}(\phi) - \widehat{T}(\phi_{0}) \mid \rightarrow 0$, for T extremal.
Note that T extremal implies $e^{\phi}/e^{\phi_{0}}=const$ (\cite{Choquet62}).

\vsp

Assume M are $C^{\infty}$ functions with monodromy, then we have $H \subset M$.
Assume C are monogenous functions, in particular $f \in C^{\infty}(x,y)$ implies $f \in C^{\infty}(x)$.
Example:  $\phi \in C^{\infty} \cap M(u,v)$ in the space of test functions for $dU \in G$ (translates) and dV conjugated in G, for instance $\phi \in \mathcal{D}_{L^{1}}$ with $H \subset \mathcal{D}_{L^{1}}$, restricted to $f \in \mathcal{D}_{L^{1}}'$. Monodromy means in particular, existence of (u,v) corresponding to G, such that $f(u,v) \in C^{\infty}$, when for instance $u \rightarrow 0$ (or $v \rightarrow 0$), that is monogenity.

\vsp

Assume $U \mathcal{L}(\phi)=\mathcal{L}(V \phi)$, that is V determines U uniquely, given a canonical topology. Assume $d V=\sigma d \omega$,  where $d \omega$ is a harmonic measure.
Conversely, when dV=0 on $\Gamma$ implies $d \omega=0$, we have a parabolic boundary.

Given monodromy, $(\phi,dN)$ can be determined from the continuation (dU,dV), that is $dN(I,0) \rightarrow d N(U,V)$ continuous, defines $(\phi,dN)$ as single valued and analytic. When $dN(I,0) \rightarrow dN(U)$ is continuous, we may have multi-valentness for $(dN,\phi)$. Conversely, $(dU,dV) \rightarrow dN$ analytic, implies that (dU,dV) can be determined locally (second order boundary condition.)

\vsp

Assume $d U \rightarrow (dU,dV) \rightarrow dN$ defines an analytic continuation, necessary for monogenity is $dN/dV \rightarrow 0$ regularly at the boundary, if further dN algebraic in dU, we have monodromy and monogenity. Assume the trajectory is given by (dU,dV). Monodromy determines which subgroup is involved. Given regularity in (u,v), sufficient for regularity in u, is a contraction domain, that is $(df/dv)(dv/du)$ $=df/du$ regular ( single valued).

When dN(U,V) is regular on a contraction domain in 3-space, we have linear continuation according to $d I \rightarrow d N$.
dN is uniquely determined by $f_{u},f_{v}$ (linearly  independent), as long as $f_{u},f_{v}$ are single valued (regular).
Example: $\mathcal{E}(dU) + \mathcal{E}(dV) \rightarrow \mathcal{E}(d U = d V) \rightarrow 0$ is exact.
Regularly situated subgroups, implies decomposability in $C^{\infty}$. Example: $\int \phi(u,u) du=\int (g_{1}(u) - g_{2}(v(u))) du$.

\subsection{A weighted condition}
\label{sec:moment}
 Assume $1/\sigma=(\vartheta + 1/\vartheta)^{-1} \in L^{1}$, with $ \sigma f  \in L^{1}$, for $f \in \mathcal{H}$.
Let $M(\vartheta)(f)=(\vartheta \tau + \vartheta^{-1} \tau^{-1})(f)$,  where $\tau$ denotes translation and $\vartheta$ is assumed real.  Assume $M(\vartheta^{-1})(f) < \sigma f(x) < M(\vartheta)(f)$
analogous with the $\rho-$ condition (\cite{Ahlfors56}), when $\vartheta$ is constant $\neq 0$.
Assume $\Sigma_{+} = \{ f(x+\eta) = f(x) \}$, that is lineality is given by $\Sigma_{+} \cup \Sigma_{-}$. Necessary for a  disk-subset, is then $\Sigma_{+}=\Sigma_{-}$. The condition above implies abscense of invariants, that is $\sigma f$ is a reduced symbol.
Example: $f(x+\eta)=f(x)=f(x-\eta)$, then $(\vartheta + 1/\vartheta) = 1 $ has bounded support and thus we do not have symmetric invariants. Note that a strict condition is necessary for inclusion, that is a majorization principle.

Consider $d M(\varphi)(f)=\varphi dU(f) + \varphi^{-1} d U^{-1}(f)$. When $dU=\rho d I$ and $dV=\vartheta dU$,
we have $d V^{-1}=(1/\vartheta \rho) d I$. Note that $\rho/\vartheta + \vartheta / \rho \in L^{1}$ does not imply $\vartheta/\rho \in L^{1}$, that is $d M(\vartheta) \in G$ does not imply $d U^{-1}VU \in G$.
Example: $d U=\rho d I$ and $dM(\vartheta^{-1})(f)$ volume preserving, that is $\vartheta^{-1}(d Uf) + \vartheta (d U^{-1} f) \sim f$, when we have $(\rho/\vartheta) dI + (\vartheta/\rho)d I \sim dI$.

\vsp

Starting from  f(x,y,z), with $\int \mid df \mid < \infty$ and a determined tangent plane,
we have a determined (single valued) normal. Note that when $V^{2}=\Omega(\varphi)$ (\cite{Lie91}), for a generatrix $\varphi$ and $dU \varphi=\Omega(\varphi)$, for instance a rotational surface, then we have on a contraction domain,
that $X dy-Y dx=0$ gives a determined tangent plane in 3-space. Further, for any function $\Phi=\Phi(\varphi)$,
$d U \Phi(\varphi)=\frac{d \Phi}{d \varphi} \Omega(\varphi)$ and $\Phi=\int \frac{d \varphi}{\Omega(\varphi)}$
(\cite{Lie91}).

Example: $d I \rightarrow (dU,dV)$  where dI is defined in the plane. The Dirichlet integral, $D_{W}$ is invariant for conformal mappings (\cite{AhlforsSario60}). Assume $\int_{W} (f_{u}^{2} + f_{v}^{2})du dv=D_{X}$,  where W=W(u,v) gives a planar representation space and X=X(x,y,z) (for instance according to Brelot). bd X limits X, this implies bd X is n-1 dimensional, in the same manner bd W is 1-dimensional. In 3-space, bd X has uniquely determined normal, given a 2-dimensional, rectifiable boundary (T(f)=$\int \mid df \mid <\infty$). For an analytic f, $D_{W}(f)=\int_{W}\mid f'(u+i v) \mid^{2} d u d v$. Example: when $T(f)=\int \mid f'(z) \mid \mid d z \mid$, over a disk $W$, we have $T(f)^{2} \leq D_{W}(f)$. However, when $D_{W}$ is defined by an extremal distribution, $<\delta(x-y),\mid f'(x) \times f'(y) \mid>$, it does not imply that T(f) has extremal definition.

\bibliographystyle{amsplain}
\bibliography{december2024}

\providecommand{\bysame}{\leavevmode\hbox to3em{\hrulefill}\thinspace}
\providecommand{\MR}{\relax\ifhmode\unskip\space\fi MR }
\providecommand{\MRhref}[2]{%
  \href{http://www.ams.org/mathscinet-getitem?mr=#1}{#2}
}
\providecommand{\href}[2]{#2}
\begin{thebibliography}{10}

\bibitem{AhlforsSario60}
L.~Sario~L. Ahlfors, \emph{Riemann Surfaces}, Princeton University Press, 1960.

\bibitem{Borel12}
E.~Borel, \emph{Les Fonctions Monog\`enes non Analytiques}, Bulletin de la
  S.M.F. \textbf{tome 40} (1912).

\bibitem{Bourbaki64}
N.~Bourbaki, \emph{Espaces Vectoriels Topologiques, livre V}, Hermann, Paris,
  1964.

\bibitem{Bourbaki89}
\bysame, \emph{General Topology, ch. 1-4}, Springer - Verlag Berlin, 1989.

\bibitem{Brelot49}
M.~Brelot, \emph{Le Probl\`eme de Dirichlet G\'eod\'esique}, C.R. (1949).

\bibitem{Chirka89}
E.~M. Chirka, \emph{Complex Analytic Sets}, Springer Science and Business
  Media, 1989.

\bibitem{Choquet62}
G.~Choquet, \emph{Le Probl\`eme des Moments.}, S\'eminaire Choquet. Initiation
  \`a l'analyse, t.1 (1962).

\bibitem{Dahn19}
T.~Dahn, \emph{Some Remarks on Polar Sets to Sums of Squares}, ArXiv (2019).

\bibitem{Dahn23}
\bysame, \emph{On Uniformity for the Polar to Partially Hypoelliptic
  Operators}, ArXiv (2023).

\bibitem{Dahn24}
\bysame, \emph{On Very Regular Representations, in Presence of Index.}, ArXiv
  (2024).

\bibitem{Riesz56}
B.~Sz.-Nagy {F. Riesz}, \emph{Functional Analysis.}, Dover Publications Inc.,
  1956.



\bibitem{Faraut70}
J.~Faraut, \emph{Semi-groupes de Mesures Complexes et Calcul Symbolique sur les
  G\'en\'erateurs Infinit\'esimaux de Semi-groupes d'Op\'erateurs.}, Annales de
  l'institut Fourier, t.20, no. 1 (1970).

\bibitem{Garding87}
L.~G{\aa}rding, \emph{Singularities in Linear Wave Propagation}, Springer
  Lecture Notes in Mathematics, 1241, Springer Lecture Notes in Mathematics,
  1241, 1987.

\bibitem{Julia19}
G.~Julia, \emph{Sur Quelques Propri\'et\'es Nouvelles des Fonctions Enti\`eres
  ou M\'eromorphes (premier m\'emoire).}, Annales scientifiques de l'E.N.S.,
  3-e s\'erie, t. 36 (1919).

\bibitem{Ahlfors56}
A.~Beurling L.~Ahlfors, \emph{The Boundary Correspondence under Quasiconformal
  Mappings}, Acta Math., Vol. 96 (1956).

\bibitem{Lelong68}
P.~Lelong, \emph{Fonctionelles Analytiques et Fonctions Enti\`eres: (n
  variables)}, Presses de l'Universit\'e de Montr\'eal, 1968.

\bibitem{BrelotChoquet51}
G.~Choquet M.~Brelot, \emph{Espaces et Lignes de Green}, Annales de l'institut
  Fourier, t.3 (1951).

\bibitem{Malgrange66}
B.~Malgrange, \emph{Ideals of Differentiable Functions.}, Tata Institute of
  Fundamental Research, Bombay, 1966.

\bibitem{Parreau51}
A.~Parreau, \emph{Sur les Moyennes des Fonctions Harmoniques et Analytiques et
  la Classification des Surfaces de Riemann.}, Annales de l'institut Fourier.
  \textbf{t.3} (1951).

\bibitem{Garding47}
L.~G{\aa}rding, \emph{Note on Continous Representations of Lie Groups}, Proc.
  Natl. Acad. Sci., USA \textbf{no. 33(11)} (1947).

\bibitem{Garding60}
\bysame, \emph{Vecteurs Analytiques dans les Repr{\'e}sentations des Groupes de
  Lie}, Bulletin de la S.M.F \textbf{t.88} (1960).

\bibitem{Oka60}
K.~Oka, \emph{Sur les Fonctions Analytiques de Plusieurs Variables.}, Hermann
  (1960).

\bibitem{Riesz15}
Hardy~G.H. Riesz, M., \emph{The General Theory of Dirichlet's Series},
  Cambridge Tracts in Mathematics and Mathematical Physics, No. 18, 1915.

\bibitem{Riesz27}
M.~Riesz, \emph{Sur les Maxima des Formes Bilin\'eaires et sur les
  Fonctionelles Lin\'eaires.}, Acta Math. 49 (1927).


\bibitem{Lie91}
S.~Lie~G. Scheffers, \emph{Vorlesungen {\"U}ber Differentialgleichungen mit
  Bekannten Infinitesimalen Transformationen.}, Teubner Leipzig, 1891.

\bibitem{Schwartz52}
L.~Schwartz, \emph{Transformation de Laplace des Distributions}, Comm. S\'em.
  Math. Univ. Lund \textbf{Tome suppl.} (1952).

\bibitem{Schwartz66}
\bysame, \emph{Th{\'e}orie des Distributions.}, Hermann, 1966.

\bibitem{Treves67}
F.~Treves, \emph{Topological Vector Spaces, Distributions and Kernels.},
  Academic Press, 1967.


\end{thebibliography}

\end{document}